\title{On the hardness of sampling independent sets beyond the tree threshold}
\author{Elchanan Mossel\thanks{Statistics, U.C. Berkeley. E-mail: {\tt
      mossel@stat.berkeley.edu}. Supported by
    an Alfred Sloan fellowship in Mathematics and by NSF grants 
    DMS-0528488, DMS-0504245 and DMS-0548249 (CAREER). Most of this
    work was done while the author was visiting Microsoft Research} \and 
    Dror Weitz\thanks {DIMACS Center, Rutgers university, Piscataway, NJ 08854
    U.S.A.\ \ Email: {\tt dror@dimacs.rutgers.edu}. Some of this work was done
    while the author was visiting Microsoft Research} \and  Nicholas Wormald 
    \thanks{Combinatorics and Optimization, University of Waterloo, Waterloo ON, 
    Canada N2L 3G1. E-mail: {\tt nwormald@uwaterloo.ca} Supported by the Canada
    Research Chairs Program and NSERC. Most of this
    work was done while the author was visiting Microsoft Research}}

\documentclass[11pt]{article}
\usepackage{amssymb,amsmath,amscd,epsfig,enumerate,graphicx}

    \newcommand{\lab}[1]{\label{#1}}           

\def\fullpage
{
\addtolength{\topmargin}{-1.5 cm}
\addtolength{\oddsidemargin}{-1.5 cm}
\addtolength{\textwidth}{+3 cm}
\addtolength{\textheight}{+3 cm}
}
\fullpage

\def\eps{\epsilon}
\newcommand{\be}{\begin{equation}}
\newcommand{\ee}{\end{equation}}
\newcommand{\bea}{\begin{eqnarray}}
\newcommand{\eea}{\end{eqnarray}}
\newcommand{\non}{\nonumber}
\newcommand{\bean}{\begin{eqnarray*}}
\newcommand{\eean}{\end{eqnarray*}}
\newcommand\eqn[1]{(\ref{#1})}
\newcommand{\bel}[1]{\be\lab{#1}}

\def\figone{

\begin{figure}[htb]
     \begin{center}
       \includegraphics[width=60mm]{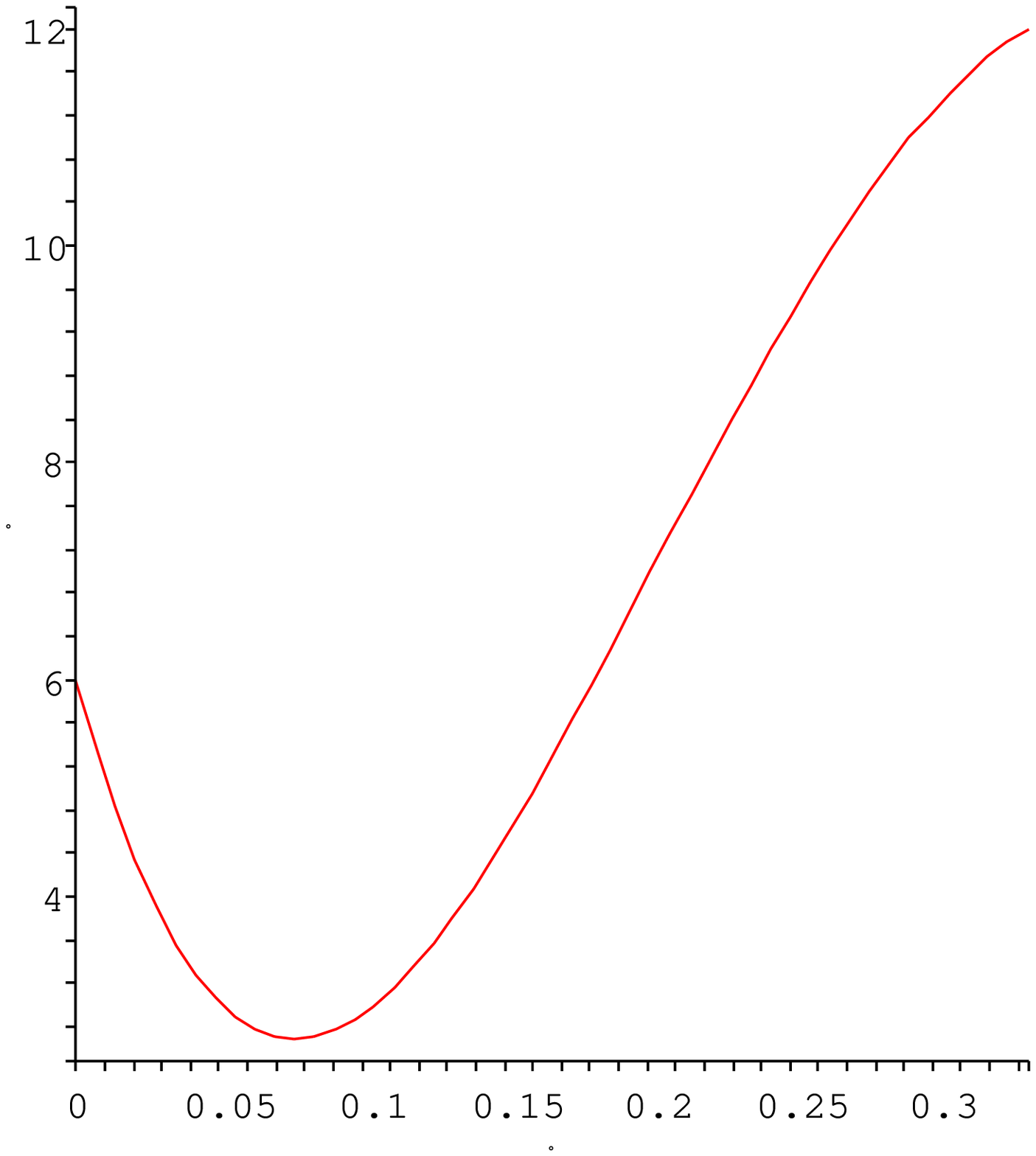}
       \hspace{2cm}
 \includegraphics[width=60mm]{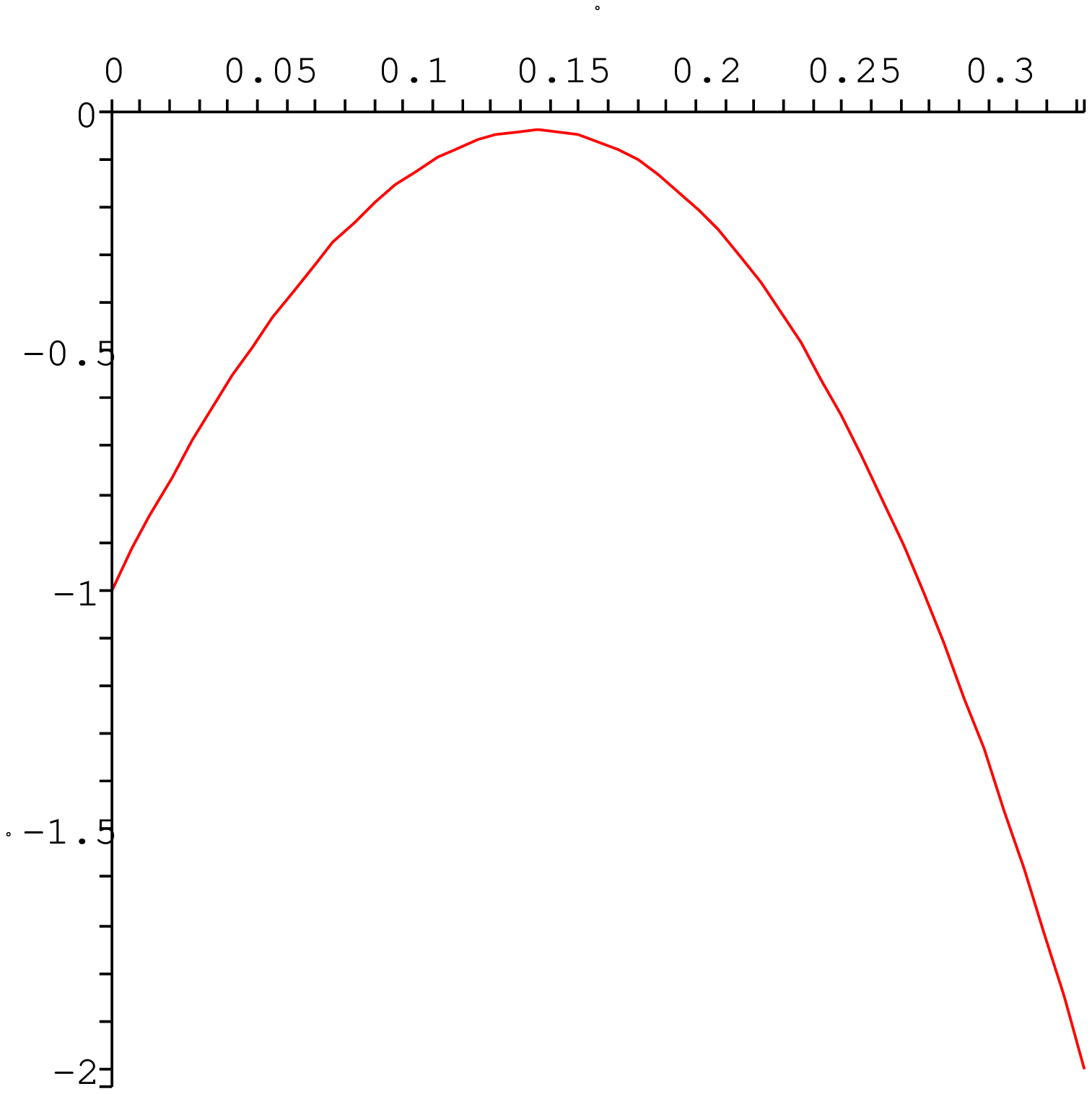}

  \caption{The two quartic polynomials}
\lab{fig1}     \end{center}

\end{figure}
}
\def\figtwo{

\begin{figure}[htb]
     \begin{center}
      \includegraphics[width=60mm]{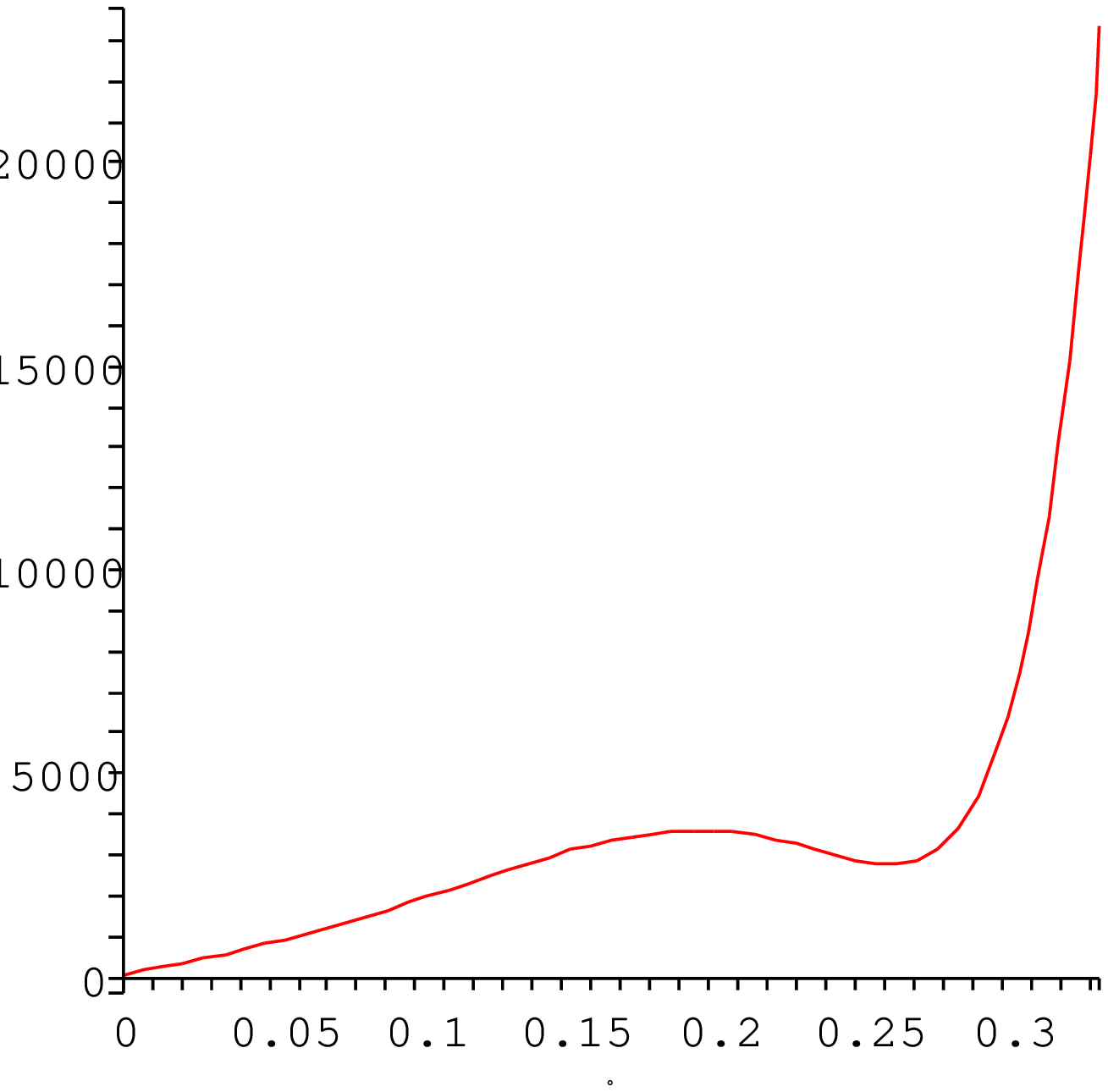}
      \hspace{2cm}
         \includegraphics[width=60mm]{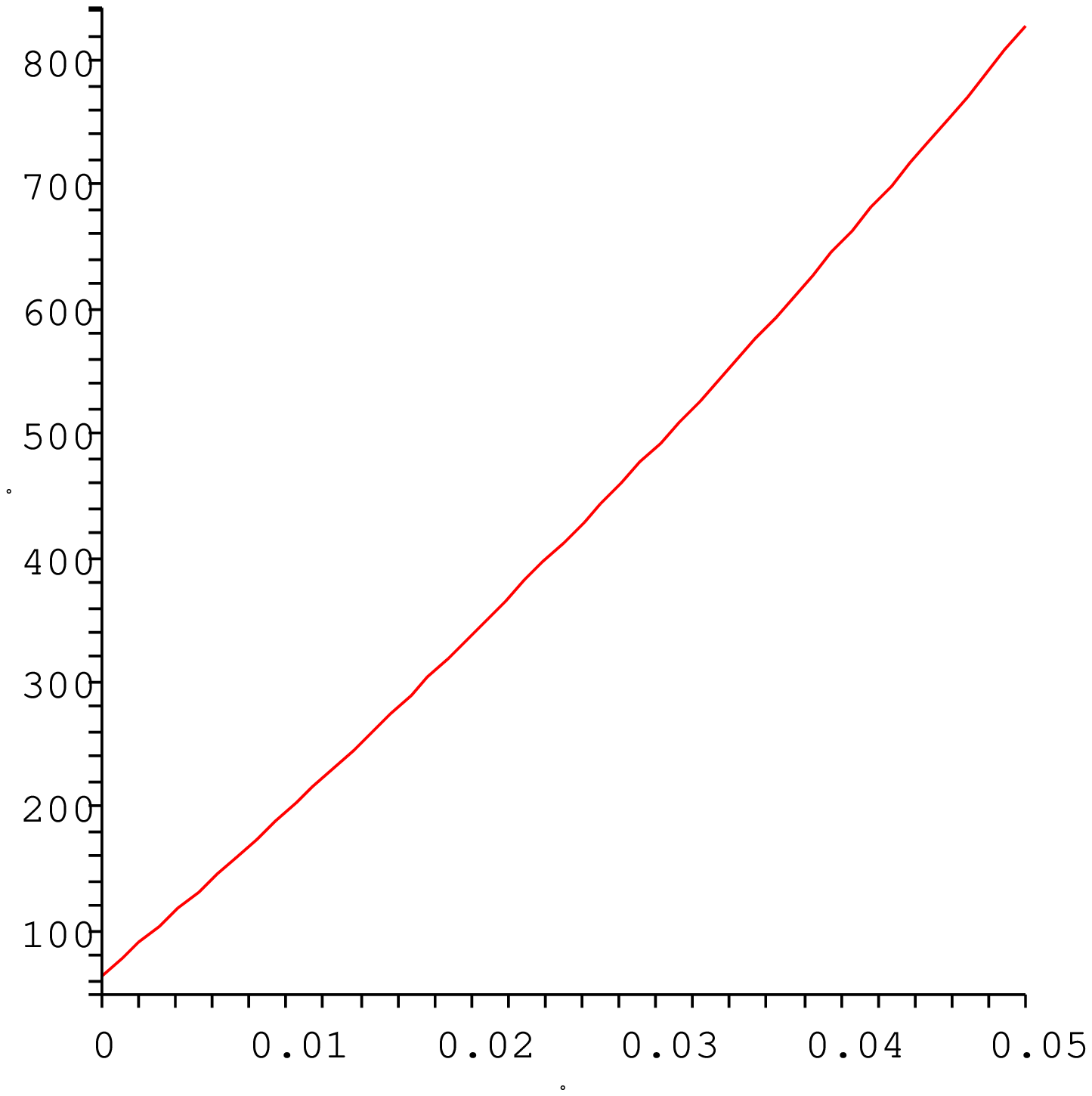}
   \caption{The degree 11 polynomial $P(\delta)$}
\lab{fig2}      \end{center}

\end{figure}
}

 \def\figthree{

\begin{figure}[htb]
     \begin{center}
      \includegraphics[width=60mm]{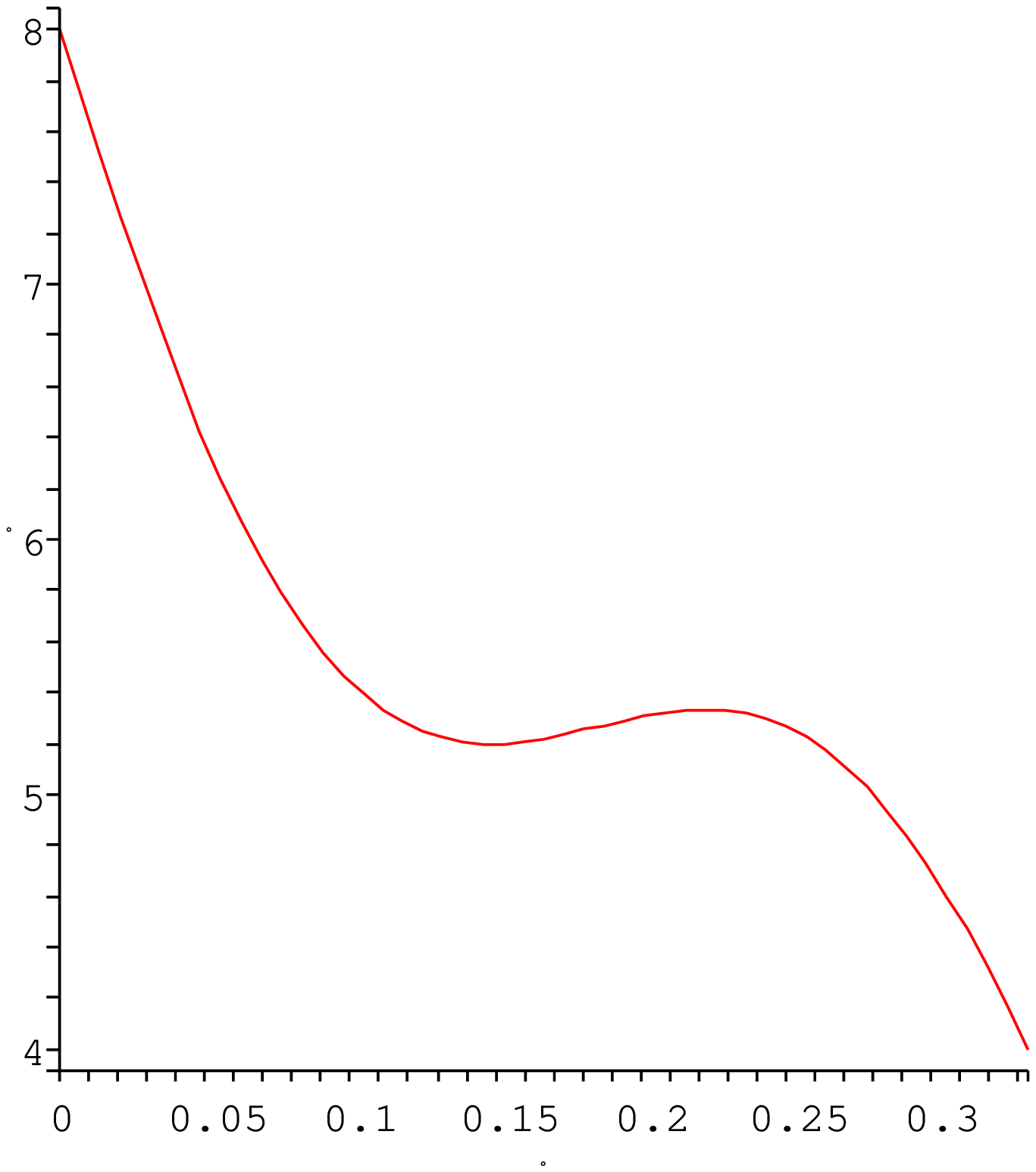}
   \caption{The quintic polynomial}

\lab{fig3}     \end{center}

\end{figure}
}

\makeatletter
\def\bigpar{\bigbreak\@afterindentfalse\@afterheading\ignorespaces}
\def\medpar{\medbreak\@afterindentfalse\@afterheading\ignorespaces}
\def\smallpar{\smallbreak\@afterindentfalse\@afterheading\ignorespaces}

\newenvironment{romanenumerate}
      {\begin{enumerate}}
      {\end{enumerate}}

\newlength{\saveindent}
\newenvironment{proof}%
      {\bigpar{\bf Proof:}\ 
             \setlength{\saveindent}{\parindent} 
                       \ignorespaces}%
      {\stopproof\ignorespaces\bigbreak \setlength{\parindent}{\saveindent}}

      {\bigpar{\bf Proof:}%
             \setlength{\saveindent}{\parindent} 
                       \,\ignorespaces}%
      {\ignorespaces\bigbreak \setlength{\parindent}{\saveindent}}

      {\bigpar{\bf Proof:}\ %
             \setlength{\saveindent}{\parindent} 
                       \ignorespaces}%
      {\ignorespaces\bigbreak \setlength{\parindent}{\saveindent}}

\newenvironment{proofof}[1]%
      {\bigpar{\bf Proof of #1:}\ %
             \setlength{\saveindent}{\parindent} 
                       \ignorespaces}%
      {\stopproof\ignorespaces\bigbreak \setlength{\parindent}{\saveindent}}

      {\bigpar{\bf Proof of #1:}\ %
             \setlength{\saveindent}{\parindent} 
                       \ignorespaces}%
      {\ignorespaces\bigbreak \setlength{\parindent}{\saveindent}}

\newenvironment{remark}%
      {\medpar{\bf Remark:}\ 
                       \ignorespaces\small}%
      {\stopproof\ignorespaces\medbreak \setlength{\parindent}{\saveindent}} 

\newenvironment{remark*}%
      {\medpar{\bf Remark:}\ 
                       \ignorespaces\small}%
      {\ignorespaces\medbreak \setlength{\parindent}{\saveindent}} 

      {\medpar{\bf Remarks:}\ 
                       \ignorespaces\small}%
      {\stopproof\ignorespaces\medbreak \setlength{\parindent}{\saveindent}}

\newenvironment{remarks*}%
      {\medpar{\bf Remarks:}\ 
                       \ignorespaces\small}%
      {\ignorespaces\medbreak \setlength{\parindent}{\saveindent}}

      {\medpar{\bf Remark:}\ %
                       \ignorespaces}%
      {\stopproof\ignorespaces\medbreak \setlength{\parindent}{\saveindent}}

      {\medpar{\bf Remarks:}\ %
                       \ignorespaces}%
      {\stopproof\ignorespaces\medbreak \setlength{\parindent}{\saveindent}}
\makeatother

\newtheorem{theorem}{Theorem}[section]   

\newtheorem{lemma}[theorem]{Lemma}
\newtheorem{claim}[theorem]{Claim}
\newtheorem{proposition}[theorem]{Proposition}
\newtheorem{corollary}[theorem]{Corollary}
\newtheorem{definition}[theorem]{Definition}

\newtheorem{conjecture}[theorem]{Conjecture}

\newtheorem{fact}[theorem]{Fact}
\newtheorem{example}{Example}[section]

      {\bigpar{\bf #1}\ %
             \setlength{\saveindent}{\parindent} 
                       \ignorespaces\it}%
      {\ignorespaces\bigbreak \setlength{\parindent}{\saveindent}}
      {\bigpar{\bf Theorem #1${}'$}\ %
             \setlength{\saveindent}{\parindent} 
                       \ignorespaces\it}%
      {\ignorespaces\bigbreak \setlength{\parindent}{\saveindent}}
      {\bigpar{\bf Theorem #1${}''$}\ %
             \setlength{\saveindent}{\parindent} \setlength{\parindent}{0pt}%
                       \ignorespaces\it}%
      {\ignorespaces\bigbreak \setlength{\parindent}{\saveindent}}

\def\begex{\begin{example}\parindent=0pt \rm}
\def\endex{\end{example}}
\def\square{\vbox{\hrule height.2pt\hbox{\vrule width.2pt height5pt \kern5pt
                                   \vrule width.2pt} \hrule height.2pt}}
\def\stopproof{\qquad\square}

\newskip\storeadskip
\newskip\storebdskip

\let\a=\alpha \let\b=\beta   \let\d=\delta  \let\e=\varepsilon
 \let\g=\gamma       \let\l=\lambda
          \let\p=\pi

\let\D=\Delta

\newcommand{\cG}{\ensuremath{\mathcal G}}

\newcommand{\cI}{\ensuremath{\mathcal I}}

\newcommand{\cM}{\ensuremath{\mathcal M}}

\newcommand{\cT}{\ensuremath{\mathcal T}}

\let\neper=e
\let\ii=i

\let\<=\langle
\let\>=\rangle

\def\half{{\textstyle{1\over2}}}

\def\Bethe{{\widehat{\Bbb T}}}

\def\pr{{\rm P}}

\newcommand{\set}[1]{\left\{#1\right\}}

\newcommand{\Var}[0]{{\mathrm{Var}}}
\newcommand{\E}[0]{\mathrm{E}}
\newcommand{\Ent}[0]{{\mathrm{H}}}
\newcommand{\RG}[0]{{\mathcal{G}}}
\newcommand{\Otilde}[0]{{\widetilde{O}}}
\newcommand{\Omegatilde}[0]{{\widetilde{\Omega}}}
\newcommand{\Thetatilde}[0]{{\widetilde{\Theta}}}
\newcommand{\pcond}[0]{{\hat{p}}}

\begin{document}

\maketitle
\date{}
\begin{abstract}
  We consider local Markov chain Monte-Carlo algorithms for sampling
  from the weighted distribution of independent sets
  with activity~$\l$, where the weight of an independent set~$I$ is
  $\l^{|I|}$. 
  A recent result has established that Gibbs sampling is rapidly mixing
  in sampling the distribution for graphs of maximum degree~$d$ and 
  $\l<\l_c(d)$, where $\l_c(d)$ is the critical activity for 
  uniqueness of the Gibbs measure
  (i.e., for decay of correlations with distance in the weighted distribution
  over independent sets) on the $d$-regular infinite tree. 

  We show that for $d \geq 3$, $\l$ just above $\l_c(d)$ 
  with high probability over $d$-regular bipartite graphs, 
  any local Markov chain Monte-Carlo algorithm
  takes exponential time before getting close to the stationary
  distribution. 

  Our results provide a rigorous justification for ``replica''  method 
  heuristics. These heuristics were invented in theoretical physics 
  and are used in order to derive predictions on Gibbs measures on 
  random graphs in terms of Gibbs measures on trees. 
  A major theoretical challenge in recent years is to provide rigorous proofs for the correctness of such predictions. Our results establish such rigorous proofs for the case of hard-core model on bipartite graphs.

  We conjecture that
  $\l_c$ is in fact the exact threshold for this computational problem, i.e.,
  that for $\l>\l_c$ it is NP-hard to approximate the above weighted sum over
  independent sets to within a factor polynomial in the size of the
  graph. 


\end{abstract}


\section{Introduction}

\subsection{Sampling weighted independent sets}
Approximately counting (or sampling) weighted independent sets is a central
problem in computational complexity, in statistical physics, where lattice gases are modelled, and in communication networks.  Typically the weights of the
independent sets are governed by an activity parameter~$\lambda$ so that the
weight of an independent set~$I$ is proportional to $\lambda^{|I|}$.

Intuitively, counting or sampling independent sets becomes harder as $\lambda$
increases. Indeed, if we could sample an independent set from the distribution
resulting from a large enough setting of~$\lambda$ then we would have an
algorithm for finding the maximum independent set of a graph, an NP-hard
problem.

In this paper we restrict our attention to local Markov chain Monte
Carlo methods for sampling weighted independent sets. 
A recent paper of the second author~\cite{W06} established the that
Gibbs sampling, a simple local Markov chain Monte Carlo algorithm is
rapidly mixing for graphs of maximum degree $d$ if 
$\l<\l_c(d)$ where   $\l_c(d)= (d-1)^{d-1}/(d-2)^d $ is the threshold 
for decay of correlations on the
regular tree of degree~$d$. In this paper we show that $\l_c(d)$ is in
fact the exact threshold for rapid mixing of local dynamics, by
showing that for random bipartite $d$-regular graphs, 
with high probability, the mixing time
of local dynamics is exponential in the size of the graph if $\l >
\l_c(d)$. 

In the following subsections we state our main result and provide
motivations for studying this problem from computational
complexity, the replica method and the role of uniqueness and
extremality 

\subsection{Our results}
In our main result we show that for $d \geq 3$ and $\l$ just above $\l_c(d)$ 
for almost all bipartite $d$-regular 
graphs, any local MCMC must take exponential
time to get close to the stationary distribution. 
More formally we show the following. 
\begin{theorem}
  \label{thm:main_simple}
  For all $d \geq 3$ there exists $\eps_d > 0$ such that for $\l_c(d) < \l < \l_c(d) + \eps_d$
  there exists $a = a(\l,d)>1$ such that with high probability
  (probability tending to $1$ as $n \to \infty$) 
for a random $d$-regular bipartite graph, 
the mixing time of the Glauber dynamics on~$G$ is $\Omega(a^n)$.
Moreover, the same claim holds true for any reversible dynamics which updates at most 
$o(n)$ nodes at each iteration.
\end{theorem}

\subsection{Computational complexity of sampling weighted independent sets}

Based on the intuition that sampling independent sets is
computatioally harder for larger values of $\l$, 
it was shown~\cite{hardness} that for any $d\ge 4$, it
is NP-hard to approximate the above weighted sum over independent sets, even to
within a polynomial factor, for graphs of maximum degree~$d$ and
$\lambda>{c\over d}$, where $c$ is a (large enough) absolute constant.

On the other hand, the existence of a fully polynomial approximation scheme has
been established for an ever-improving sequence~\cite{hardness,dyer_ind,V01,W06} of
bounds on~$\l$ which are also inverse linear in the maximum degree. It is a
fascinating challenge to determine the exact threshold (in terms of $\l$ as a
function of $d$) for which the counting is hard to approximate. (It is known
that for any given $(d,\l)$, either there exists a fully polynomial
approximation scheme or it is NP-hard to approximate the sum over independent
sets to within a polynomial factor~\cite{sinclair}).
 
It has been speculated (though not formally conjectured) that the  
hardness threshold corresponds to the threshold for decay of spatial
correlations for the weighted distribution over independent sets (uniqueness of
the {\em Gibbs} measure). In particular, it is was speculated that the critical
activity $\l_c(d)= (d-1)^{d-1}/(d-2)^d$ for decay of correlations on the
regular tree of degree~$d$ is also the threshold for the computational problem.
This is supported by a recent paper of the second author~\cite{W06} 
which established the
existence of a fully polynomial approximation scheme for counting independent
sets for every~$d$ and $\l<\l_c(d)$.  The present paper is motivated by this
speculation. In particular, we conjecture the following.
\begin{conjecture}
  \label{c:hardness}
  For every $d>3$ and all $\l>\l_c(d)$, unless $RP=NP$ there does not exist a
  fully polynomial approximation scheme for counting weighted independent sets
  of graphs of maximum degree~$d$ with activity~$\l$.
\end{conjecture}
We provide evidence supporting the above conjecture by analyzing local Markov
chain Monte Carlo algorithms for sampling independent sets. By local algorithms
we refer to algorithms in which the number of vertices updated in a single step
of the chain is $o(n)$, where $n$ is the number of vertices of the graph.

Of course, slow
mixing of local MCMC algorithms as established in this paper does not
generally imply the hardness result stated in the conjecture. For
example, the {\em Ising
  model} at low temperatures is an example where the local MCMC algorithm mixes
slowly while there is an FPRAS~\cite{JS} for computing the partition function.

However, we believe that our results do give support for the
conjecture.  First, unlike the Ising model where there is no computational
phase transition and an FPRAS exits for all temperatures~\cite{JS}, we know (as
already mentioned above) that approximate counting of independent sets is computationally hard
for $\lambda\ge {c\over d}$ for a large enough absolute constant $c$.  Second,
 our result establishes that for $\l$ just above $\l_c(d)$ and
most $d$-regular bipartite graphs {\em balanced} independent sets form an
exponentially small bottleneck, i.e., the density of a typical independent set
in these graphs has a positive bias to either $V_1$ or $V_2$ except with
exponentially small probability. The availability of such graphs
could open the way for hardness constructions in which these graphs would be
used as gadgets for, e.g., encoding a binary variable, and thus that the
hardness threshold for approximate counting of independent sets coincides with
the threshold for the availability of such graphs. See~\cite{DFJ} for
an easier 
construction of a similar flavor that was used to establish that approximate
counting of independent sets with $\l=1$ is hard for $d\ge 24$.

\subsection{Our result and the Replica method}
Another important motivation for our result is the ``replica'' 
heuristic developed in theoretical 
physics~\cite{MezardParisi:86,MezardParisi:87}. 
This method gives predictions on the behavior of Gibbs measures on random 
graphs that are based on the analysis of Gibbs measures on trees. 
The method has been extensively used to derive predictions of the behavior of 
many random systems~\cite{MePaVi:91}. The method has also yielded an 
empirically  effective algorithm for solving random 3-SAT 
problems at the highest known densities~\cite{MePaZe:03}. 

The theoretical study of the replica method has been a major challenge in 
mathematics, theoretical physics, probability, engineering and computer 
science. A number of results proved the validity of the method for various 
specific models such as the SK model~\cite{Talagrand:03}, 
the assignment problem~\cite{Aldous:01} 
and some results in coding theory~\cite{Richardson01b,Lubyetal}. 
However, all of these results deal with specific problems. No general results 
are known for the applicability of the method.

Our results provide another example where it can be rigorously shown that 
replica calculations do determine the behavior of Gibbs measures on 
random graphs. In particular, our results are strong enough to actually determine the dynamics properties of the Gibbs measures, i.e., the convergence time of reversible local dynamics. 

\subsection{Uniqueness, Extremality and Slow Mixing}
Another important aspect of our results is that establishing exponentially slow mixing
for random $d$-regular graphs for $\l$ just above $\l_c$ shows that these graphs,
while locally similar to the $d$-regular tree, behave in a very different manner
than the tree w.r.t.\ the mixing time of the chain for $\l>\l_c$.  On the
regular tree, the mixing time of the local Markov chain is $O(n\log n)$ even
above $\l_c$. More precisely, there exists a second threshold~\cite{MSW04}
$\l_1>\l_c$ such that the mixing time is $O(n\log n)$ for $\l\le \l_1$. (On the
tree the mixing time is polynomial in~$n$ for any $\l$~\cite{BKMP05}.)  This
corresponds to the fact that the Gibbs measure is {\em extremal} for $\l\le
\l_c$ (see also~\cite{winkler}).

Our results show that while random regular graphs are locally
tree-like, they have very different properties when it comes to
convergence of local dynamics. In particular, while the threshold for
uniqueness of the Gibbs measure (appropriately defined) 
is the same in random graphs and the tree, in
random graphs --- unlike the tree --- this threshold is also where the Gibbs
measure ceases to be extremal and where the mixing time undergoes a sharp
transition and becomes exponentially slow. This is in line with the intuition
that the mixing time depends more crucially on the size of the
separators of the graph than on their local structure. 
Trees have very small separators (e.g., the root) and thus have fast 
mixing time even above $\l_c$.  
Random graphs on the other hand are expanders, and indeed
have exponentially slow mixing already just above $\l_c$.

\subsection{Proof Technique}
Our proof borrows its initial approach from the proof in~\cite{DFJ} where the case $\l=1$ is
analyzed. It is shown in~\cite{DFJ} that most random $d$-regular bipartite
graphs exhibit an exponentially small bottleneck for $d\ge 6$, causing the
Markov chain to take exponential time to mix.  Note that $\l_c(d)\le \l_c(6)
\approx 0.763 <1$ for every $d\ge 6$ so even though the latter result is tight
for $\l=1$ in terms of the degree $d$, it is far from tight for general $\l$ as
a function of $d$.

In this paper we show that the model suggested in~\cite{DFJ} can be in fact be
analyzed all the way down to $\l_c(d)$ for $d \geq 3$.  
We note that the arguments
of~\cite{DFJ} do not extend to give a proof of slow mixing all the way down to
$\l_c(d)$: 
in that case, where a lower bound on the number of balanced independent sets of a given size and location was required, it sufficed to use a crude lower bound that applied for all graphs. This was made easier because the balanced independent sets were far from the typical ones (in terms of their ``balancedness"). For our present purposes, we need to consider cases where the typical independent sets are quite close to being balanced. As a result, we found it necessary to use some sophisticated second moment calculations. This is even just to show the existence of at least one graph with an exponentially small bottleneck of the type we desire. To show that almost all graphs have such a bottleneck we used the small subgraph conditioning method of the third author and Robinson (see~\cite{models}).

These apply for all $\l$ slightly above $\l_c(d)$. 
We believe that the same result holds for
for all $d$ and $\l>\l_c(d)$. Also, we would like to point out 
that our first-moment
analysis, while similar to~\cite{DFJ}, 
exhibits more explicitly the role of Gibbs measures on the
regular tree.

\noindent
{\bf Acknowledgments:} E.M. and D.W. wish to thank Alistair Sinclair
for interesting discussions. 

\section{Preliminaries and statements of results}\lab{s:pre}

\subsection{The random graph model}
We consider the following model for random graphs. The graphs are all bipartite
with vertex sets $V_1,V_2$ of size~$n$ each. We choose $d$~random perfect
matchings between the two vertex sets so that every vertex has degree~$d$. Note
that there maybe parallel edges (with probability asymptotic to a constant less than 1).  However, for the
hard core model, a parallel edge has exactly the same effect as does a single
edge. So this can alternatively be regarded as a model of simple graphs in which
the degree of each node is at most $d$.  We denote the above distribution over
graphs by $\RG\equiv\RG(n,d)$. 
As is common in the theory of random
graphs, we will use the term {\em asymptotically almost surely (a.a.s.)} 
to refer to a sqequence of
probabilities converging to $1$ as $n \to \infty$. 

\begin{remark}
The probability space $\RG(n,d)$ we are working with is the set of
bipartite (multi)graphs obtained by taking $d$ random perfect
matchings between two sets $V_1$ and $V_2$ of $n$ vertices each. 
This probability space is contiguous with a uniformly random
$d$-regular graph 
(see the note after the proof of~\cite[Theorem 4]{MRRW}),  
and hence, all of the results below proven a.a.s.\ for $\RG(n,d)$ 
consequently also hold
a.a.s.\ in the uniform model, with the uniform distribution 
over simple bipartite $d$-regular
graphs  
(as well as various other models contiguous to it).
\end{remark}

\subsection{The hard-core model}
Let $\lambda>0$ be an activity parameter and~$G$ a finite graph. Denote
by~$\cI_G$ the set of independent sets of~$G$. The hard-core
distribution on~$G$ with activity $\lambda$, denoted $\mu_{G,\lambda}$, is the
distribution over $\cI_G$ in which the probability of an
independent set $I\in\cI_G$ is proportional to $\lambda^{|I|}$, i.e.,
\begin{equation}
  \label{eq:Gibbs}
  \mu_{G,\lambda}[I]={\lambda^{|I|}\over Z_{G,\lambda}},
\end{equation}
where $Z_{G,\lambda} = \sum_{I\in\cI_G} \lambda^{|I|}$.

\subsection{Gibbs measures on the regular tree}
Let $\Bethe^d$ be the infinite regular tree of degree~$d$. A probability
measure~$\mu$ on independent sets of~$\Bethe^d$ is Gibbs if for every finite
subtree~$T$, conditioning~$\mu$ upon the event that all the vertices on the outer
boundary of~$T$ are unoccupied gives the same probability distribution on
independent sets of~$T$ as defined by~(\ref{eq:Gibbs}) with $G=T$. Moreover, $\mu$ is a {\em
  simple} Gibbs measure on~$\Bethe^d$ if for any vertex~$v$, conditioning~$\mu$
on any of the two possible values at~$v$ results in a measure in which the
configurations on the $d$ (infinite) subtrees rooted at the children of~$v$ are
independent of each other. (Notice that the probability distribution on a finite
subtree is always simple.) A {\em translation-invariant} Gibbs measure
on~$\Bethe^d$ is a measure that is invariant under all automorphisms
of~$\Bethe^d$. Finally, a {\em semi}-translation-invariant Gibbs measure
on~$\Bethe^b$ is one that is invariant under all parity-preserving automorphisms
of~$\Bethe^d$.

It is well known~\cite{kelly} that the hard-core model on~$\Bethe^d$ admits a
unique simple translation invariant Gibbs measure for all values of~$\l$. For
$\l\le \l_c(d)\equiv  (d-1)^{d-1}/ (d-2)^d $ this measure is the unique
Gibbs measure of any kind. However, for $\l>\l_c$ there are additional Gibbs
measures, and in particular, two additional simple semi-invariant measures in
which the vertices of one parity are more likely to be occupied than the
vertices of the other parity.

\subsection{The Glauber dynamics and other local dynamics}
Even though our results apply to any Markov chain Monte-Carlo algorithm that
updates at most $o(n)$ vertices in one step, for convenience and definiteness,
we will first discuss a well-known and simple Markov chain for
sampling weighted independent sets, called the {\em Glauber dynamics}. This
chain is defined as follows. Starting from the current independent set~$I$,
choose u.a.r.\ a vertex~$v$ from $V_1\cup V_2$. If all the neighbors of~$v$ are
unoccupied, set~$v$ to be occupied with probability $ \l/( \l+1)$ and
otherwise set~$v$ to be unoccupied. It is easy to verify that this Markov chain
converges to the hard-core distribution $\mu_{G,\l}$. 

More generally, we will consider dynamics where at each stage a
(random) set of vertices $W$ is chosen according to some fixed
distribution. Then the configuration of the vertices at $W$ is sampled
according to the conditional probability at $(V_1 \cup V_2) \setminus
W$. We will only consider cases where the sets $W$ chosen are of size
$o(n)$. It is easy to see under mild conditions that this dynamics
also converges to the hard-core distribution $\mu_{G,\l}$.  

The main question for both dynamics 
concerns how many steps it takes the chain to get sufficiently close to
this stationary distribution.  The {\em mixing time }of the chain is defined
as the number of steps needed in order to guarantee that the chain, starting
from an arbitrary state, is within total variation distance $1/2e$ from the
stationary distribution.

Our method for establishing slow mixing of the Glauber dynamics is based on
conductance type arguments. Namely, in order to prove slow mixing we will show
the existence of $A\subset\cI$ whose measure is at
most $\half$ and whose {\em boundary} is exponentially smaller, i.e., the
probability of escaping~$A$ is exponentially small. The existence of such a
subset is well-known to imply slow mixing of the Markov chain. For example, the
following is taken from~\cite[Claim~2.3]{DFJ}.
\begin{claim}
  \label{cl:slow_mixing_via_bottleneck}
  Let $\cM$ be a Markov chain with state space $\Omega$, transition matrix~$P$,
  and stationary distribution~$\mu$. Let $A\subset\Omega$ be a set of states
  such that $\mu[A]\le \half$, and $B\subset \Omega$ be a set of states that form
  a ``barrier'' in the sense $P_{ij}=0$ whenever $i\in A\setminus B$ and $j\in
  A^c\setminus B$. Then the mixing time of $\cM$ is at least ${\mu[A]\over
    8\mu[B]}$.
\end{claim}

\subsection{Main result}
We show that for $d\ge 3$ and $\l$ just above $\lambda_c$, 
a.a.s.\ the mixing time of a random graph  drawn from $\RG(n,d)$ is 
exponential in~$n$. As explained above, this is done by establishing an
exponentially small bottleneck in the state space. For a given graph~$G$, let
$\cI_1 = \set {I\in \cI_G \,|\, |I\cap V_1| > |I\cap V_2|}$, define $\cI_2$
similarly, and let $\cI_B = \set {I\in \cI_G \,|\, |I\cap V_1| = |I\cap V_2|}$.
Notice that $\cI_B$ forms a barrier between $\cI_1$ and $\cI_2$, i.e., the
Markov chain must go through $\cI_B$ in order to cross from $\cI_1$ to $\cI_2$
and vice versa. 

Similarly, given $\tau > 0$ we define 
$\cI_1^{\tau} = \set {I\in \cI_G \,|\, |I\cap V_1| > |I\cap V_2| +
  \tau n}$, define $\cI_2^{\tau}$
similarly, and let 
$\cI_B^{\tau} = \set {I\in \cI_G \,|\, |I\cap V_1| - |I\cap V_2| \leq
  \tau n}$.
Now $\cI_B^{\tau}$ forms a barrier between $\cI_1^{\tau}$ and
$\cI_2^{\tau}$ for local Markov chains. In other words, for any reversible
chain that updates at most $\tau n$ vertices at an iteration, the
Markov chain must go through $\cI_B^{\tau}$ in order to cross from 
$\cI_1^{\tau}$ to $\cI_2^{\tau}$
and vice versa.

Our main result establishes that for a random graph from $\RG$,
$\cI_B$ is an exponentially small bottleneck.
\begin{theorem}
  \label{t:bottleneck}
  There exists function $a(d)$ and $\eps_d$ defined for $d \geq 3$ and 
  satisfying 
  \begin{itemize}
  \item
  $a(d) > 1$ for all $d \geq 3$.
  \item
  $\eps_d > 0$ for all $d \geq 3$.
  \end{itemize}
  such that for all $\l_c(d) < \l < \l_c(d) + \eps_d$ 
  there exists $a = a(\l,d) >1, \delta = \delta(\l,d) > 0$ 
  such that a.a.s., a graph 
  drawn from $\RG(n,d)$ satisfies 
  $\mu[\cI_B] \le a^{-n} \min\set{\mu[\cI_1],\mu[\cI_2]}$. Moreover
  $\mu[\cI_B^{\delta}] \le a^{-n} \min\set{\mu[\cI_1^{\delta}],\mu[\cI_2^{\delta}]}$.
   
\end{theorem}
Since $\min\set{\mu[\cI_1],\mu[\cI_2]} \le \half$,
applying Claim~\ref{cl:slow_mixing_via_bottleneck} gives the following. 
\begin{corollary}
  \label{cr:main_d=3}
  For all $d \geq 3$ there exists $\eps_d > 0$ such that for $\l_c(d) < \l < \l_c(d) + \eps_d$
  there exists $a = a(\l)>1$ such that a.a.s. 
 for a graph $G$ drawn from $\RG(n,d)$, the mixing time of the Glauber dynamics
  on~$G$ is $\Omega(a^n)$.
The same claim holds true for any reversible dynamics which updates at most 
  $o(n)$ nodes at each iteration.
\end{corollary}
In order to establish the above result we needed to resort to certain detailed
calculations which we only carried out for $\l$ close to $\l_c$. However,
we believe the same result holds for every $d\ge 3$ and any $\l>\l_c(d)$.
\begin{conjecture}
  \label{c:bottleneck}
  For every $d\ge 3$ and any $\l>\l_c(d)$  
  there exists $a = a(d,\l) >1$ such that a.a.s. on   
  $\RG(n,d)$ it holds that $\mu[\cI_B] \le a^{-n}
  \min\set{\mu[\cI_1],\mu[\cI_2]}$.
\end{conjecture}

\section{Proof of the main theorem}
In this section we prove Theorem~\ref{t:bottleneck}, which states that for $\l>\l_c$,
the set of balanced independent sets forms an exponentially small bottleneck.
The analysis proceeds by calculating for any given pair of densities
$(\alpha,\beta)$ the weight of independent sets that occupy $\alpha n$ and
$\beta n$ vertices of $V_1$ and $V_2$, respectively.  Naturally we 
assume that $\a n$ and $\b n$ are both integers. 
Let $\cI_G^{\alpha,\beta} = \set
{I\in\cI \,|\, I\cap V_1 = \a n, I\cap V_2=\b n}$ and $Z_G^{\a,\b} =
\sum_{I\in\cI_G^{\a,\b}} \lambda^{(\a+\b)n}$.  Below we use $\Otilde(\cdot)$ and
$\Omegatilde(\cdot)$ to denote $O(\cdot)$ and $\Omega(\cdot)$ up to factors that
are polynomial in~$n$.

Our first step is analyzing the expected weight (over graphs) of the different
possible densities. In what
follows, $\Ent(x)$ is the entropy function w.r.t.\ natural logarithms, i.e.,
$\Ent(x)=-x\ln(x)-(1-x)\ln(1-x)$. The following proposition is proved
in Section~\ref{sec:first}.
\begin{proposition}
  \label{p:first_mom}
  Fix $d$ and $\l$. For any $\a,\b\in[0,1]$,
  $$\E_\RG[Z_G^{\a,\b}] \;=\; \Thetatilde(e^{\Phi_1(\a,\b)n}),$$
  where 
  \begin{equation*}
    \label{eq:first_mom}
    \Phi_1(\alpha,\beta) \;=\; (\alpha+\beta)\ln(\lambda) + \Ent(\alpha)+\Ent(\beta) +
    d\cdot \Psi_1(\a,\b),
  \end{equation*}
  and
  $$
  \Psi_1 \;=\; (1-\beta)\Ent({\alpha\over 1-\beta}) - \Ent(\alpha).
  $$
\end{proposition}

The analytic properties of $\Phi_1$ in the triangle
$$\cT = \set {(\a,\b) \,|\, \a,\b\ge 0 \mbox{ and } \a+\b\le 1}$$
play an important role similarly to~\cite{DFJ} where it was needed for the case $\l=1$. 

The analysis for general
$\l$ is essentially the same. However, our analysis demonstrates explicitly the 
role of  
Gibbs measures for the hard-core model on the infinite $d$-regular
tree. In Section~\ref{sec:first} we derive the following lemma. 
\begin{lemma} \label{lem:phi}
  \label{l:phi1}
  The following holds:
  \begin{itemize}
  \item
    The maximum of $\Phi_1$ over $\cT$ along the line $\a=\b$ is achieved at
    $\a=\b=p^*$, where $p^*$ is the probability that any given vertex is occupied
    in the unique simple translation-invariant Gibbs measure on
    $\Bethe^d$.  
  \item
    If $\l\le\l_c(d)$ then $\a=\b=p^*$ is also the unique maximum of $\Phi_1$ over
    the whole of $\cT$. 
  \item
    If $\l>\l_c(d)$ then $\a=\b=p^*$ is a saddle point and
    there exist $p_1<p^*<p_2$ for which the only two maxima of $\Phi_1$ over $\cT$
    are at $(\a=p_1,\b=p_2)$ and $(\a=p_2,\b=p_1)$. Furthermore, $p_1,p_2$
    are the probabilities of occupancy of even and odd vertices, respectively, in
    the two additional simple semi-translation-invariant Gibbs measures
    on~$\Bethe^d$. 
  \item 
    $p^{\ast},p_1$ and $p_2$ all vary continuously with $\l$. 
    Moreover, $p_1 - p^{*},p_2 - p^{*} \to 0$ 
    as $\l \to \l_c$ from above.  
  \end{itemize}
\end{lemma}

Note that if $\Phi_1(\a,\b)$ was the exponent of the weight of
$(\a,\b)$-sets in a {\em specific} graph and not only as expected value then
this graph would exhibit the bottleneck as stated in Theorem~\ref{t:bottleneck}.
Indeed, suppose that for $\l>\l_c(d)$ there exists a graph $G$ with
$|V_1|=|V_2|=n$ and maximum degree~$d$ for which
$Z_G^{\a,\b}=\Thetatilde(e^{\Phi_1(\a,\b)n})$. It would then follow that
${\mu[\cI_B]}/{\mu[\cI_1]} \le {n Z^{p^*,p^*}}/{ Z^{p_2,p_1}} = n\cdot
\Thetatilde (e^{n[\Phi_1(p^*,p^*)-\Phi_2(p_2,p_1)]}) = \Thetatilde(a^{-n})$, where
$a=e^{\Phi_1(p^*,p^*)-\Phi_2(p_2,p_1)} > 1$. Similarly, it would also follow
that ${\mu[\cI_B]}/{\mu[\cI_2]} \le \Thetatilde(a^{-n})$, and the combination of
the two inequalities would establish the existence of the desired bottleneck.
  
Thus, if we could establish concentration for the random
variables $Z^{\a,\b}$, then Theorem~\ref{t:bottleneck} would follow. Note that
since these variables are exponentially large, upper bounds are easily
derived using Markov's inequality. Specifically, since the expected value of the
total weight of balanced sets is at most $n\cdot \Otilde(e^{\Phi_1(p^*,p^*)n})$
then for most graphs in $\RG(n,d)$ the total weight of the balanced sets is
$\Otilde(e^{\Phi_1(p^*,p^*)n})$. 

If we could additionally show that for $\l>\l_c(d)$ and
for most graphs~$G$ in $\RG(n,d)$, 
both $Z_G^{p_1,p_2}$ and $Z_G^{p_2,p_1}$
are $\Omegatilde(e^{\Phi_1(p_1,p_2)n})$, i.e., these variables are within
subexponential factors from their expected values, then the proof
would follow. 

It is clear that such concentration cannot hold for general
$(\a,\b)$. For example if $d$ is large and $\a=\b=1/2$ then $Z_G^{\a,\b}$
will be $0$ except with exponentially small probability while its
expected value will still be exponentially large. 

However, we conjecture that concentration does hold for some values of 
$(\a,\b)$, including the relevant values
$(p_1,p_2)$. In particular, we show that for $d \geq 3$, 
this concentration holds in the neighborhood of the point $\a=\b=1/d$.  
The proof of the required concentration follows in two stages which we
first describe somewhat roughly and then state them more formally. 
\begin{itemize}
\item First we show that for values of $\a$ and $\b$ close to $1/d$, the second
moment of $Z_G^{\a,\b}$ is of the same order as the square of the
first moment of $Z_G^{\a,\b}$. In fact, we calculate the limiting
ratio as $n \to \infty$. Using the exact calculation and the second
moment method actually allows one to obtain  weaker version of 
Theorem~\ref{t:bottleneck}, where the result is obtained with
positive probability, say probability at least $1/2$, instead of a.a.s. 
\item To obtain a high probability result we use the {\em small graph
    conditioning method}. This method ``explains" the variance of 
$Z_G^{\a,\b}$ by the interaction between the numbers of short cycles
    and the random variable $Z_G^{\a,\b}$. 
For background on the conditioning method, 
see \cite[Theorem 9.12--Remark 9.18]{JLR} and
\cite[Theorem 4.1]{models}. 
\end{itemize}

\begin{theorem} \label{thm:2nd}
 There exists an $\eta(d) > 0$ such that for fixed $\a$ and $\b$
 satisfying $|\a - 1/d| < \eta$ and $|\b - 1/d| < \eta$ it holds that 
  $$
   \lim_{n \to \infty} 
    \frac{\E_{\RG(n,d)}[(Z_G^{\a,\b})^2]}{\E^2_{\RG(n,d)}[Z_G^{\a,\b}]}
    =\tau^{\a,\b}(d)
$$
(the limit is taken over all $n$ such that $n\a,n\b$ are integer), 
where
$$
\tau^{\a,\b}(d)=
\frac{  (1-\a-\b-\a\b)^d}{(1-\a-\b+2\a\b)^{(d-1)/2} (1-\a-\b)^{(d-1)/2}
 \sqrt{1-\a-\b+d\a\b}\sqrt{1-\a-\b -(d-2)\a\b }  }\, . 
$$
\end{theorem}
Theorem~\ref{thm:2nd} is proven in Section~\ref{sec:logsec} and 
Section~\ref{sec:ratio}. Very roughly speaking, the main step of the proof is 
showing that for values of $\a,\b$ close to $1/d$ the major
contribution to the second moments comes from pairs of independent
sets of size $(\alpha,\beta)$ that are ``uncoupled''. This means for
example that the intersections of these sets are of sizes $\alpha^2$
and $\beta^2$ respectively. Moreover, for each of the $d$ matchings
defining the graph, the number of edges between the $\alpha$ size set
in one copy and a $\beta$ size set in other copy is of size $\alpha
\beta$. Once this is established, using Gaussian integration
one obtains Theorem~\ref{thm:2nd}. 

\begin{theorem} \label{thm:aas}
 There exists an $\eta(d) > 0$ such that for fixed $\a$ and $\b$
 satisfying $|\a - 1/d| < \eta$ and $|\b - 1/d| < \eta$ it holds
 a.a.s. that that 
\[
Z_G^{\a,\b} \geq \frac{1}{n} \E_{\RG(n,d)}[Z_G^{\a,\b}]
\]
(the inequality holds for values of $n$ such that $n\a,n\b$ are integer).
\end{theorem}

The proof of Theorem~\ref{thm:aas} in Section~\ref{sec:aas} 
uses the small graph conditioning
method. This method requires in particular the exact value of
$\tau^{\a,\b}(d)$. We now use Theorem~\ref{thm:aas} 
to prove the main result, Theorem~\ref{t:bottleneck}.
\begin{proofof}{Theorem~\ref{t:bottleneck}}
Since $\l > \l_c(d)$, by Lemma~\ref{lem:phi} we 
may let $(\a,\b)=(\a,\b)(\l)$ be the maximum point of $\Phi_1$ in which
$\a<\b$. Since $\Phi_1$ is continuous we may choose $\delta > 0$ such
that if 
\[
\Gamma :=
\min_{|x-\a| < \delta|, |y-\b| < \delta} \Phi_1(x,y),\,\,\, 
\Delta := 
\max_{x,y : |x-y| \leq \delta} \Phi_1(x,y) 
\]
then 
\[
\eps := \Gamma - \Delta > 0.
\]
By Markov inequality we have 
\begin{equation} \label{eq:markov}
\sum_{x,y : |x-y| \leq \delta} Z_G^{x,y} \leq \exp(n(\Delta +
\frac{\eps}{4}))
\end{equation}
(note that there are at most $n^2$ terms in the sum above). 
Note furthermore that without loss of generality we may chose $\delta$
small enough so that the two sets 
\[
\{(x,y) : |x-\a| \leq 2\delta|, |y-\b| \leq 2\delta\}, \,\,\, 
\{(x,y) : |x-y| \leq 2\delta\} 
\]
are disjoint and $\delta \leq \eta/4$ where $\eta = \eta(d)$ is used
in the statement of Theorem~\ref{thm:aas}. 

By Lemma~\ref{lem:phi} $(\a,\b)$ is a continuous function of
$\l$. Therefore there exists an $\eps_d$ such that if 
$\l_c(d) < \l < \l_c(d) + \eps_d$ then $|\a-1/d| < \eta/4$ and
$|\b-1/d| < \eta/4$. We may now chose
$(\a_1,\b_1),\ldots,(\a_k,\b_k)$ all satisfying $|\a_i - \a| < \eta/4$
and $|\b_i - \b| < \eta/4$ 
such that for all $n$ large enough we have that one of 
$n (\a_1,\b_1), \ldots, n (\a_k, \b_k)$ is an integer point. 
From Theorem~\ref{thm:aas} it follows that a.a.s.
\begin{equation} \label{eq:aas1}
\exists i, Z_G^{\a_i,\b_i} \geq \exp(n(\Gamma - \frac{\eps}{4}))
\end{equation}
and a similar statement holds for $(\b_i,\a_i)$ instead of
$(\a_i,\b_i)$. 
By~(\ref{eq:markov}) we have that a.a.s. 
\[
\mu(\cI_B^{\delta}) \leq \exp(n(\Delta + \frac{\eps}{4})) 
\]
and by~(\ref{eq:aas1}) we have that a.a.s.
\[
\min\set{\mu(\cI_1^{\delta}),\mu(\cI_2^{\delta})} \geq 
\exp(n(\Gamma - \frac{\eps}{4})).
\]
Therefore we conclude that a.a.s. 
\[
\mu[\cI_B^{\delta}] \le \exp(-\frac{\eps n}{4}) 
\min\set{\mu[\cI_1^{\delta}],\mu[\cI_2^{\delta}]}
\]
as needed.
\end{proofof}

\section{Logarithm of the first moment} \label{sec:first}
Here we prove the claims made regarding the expected 
total weight of independent sets with a given pair of densities. The
proof below emphasizes the role of the hard core model on the tree.
\begin{proofof}{Proposition~\ref{p:first_mom}}
  Notice that for a given subset $I$ that occupies $\alpha n$ and $\beta n$
  vertices of $V_1$ and $V_2$, respectively, the probability that $I$ is an
  independent set of the chosen graph~$G$ is simply the probability that all the~$d$
  chosen matchings do not connect the subset of size $\a n$ (of $V_1$) with the
  subset of size $\beta n$ (of $V_2$), i.e., in each of the matchings all the
  $\alpha n$ edges connected to the $\alpha$ subset  fall  outside the
  $\beta$ subset. Summing all subsets of density $(\alpha,\beta)$ gives:
  \begin{eqnarray*}
    \E_\RG [Z_G^{\alpha,\beta}] & = &
    {n \choose \alpha n}{n \choose \beta n}\lambda^{(\alpha+\beta)n}\left[
      {{(1-\beta)n\choose \alpha n} \over {n\choose \alpha n}}\right]^d\\
    & \approx &
    \exp \left\{n\left[(\alpha+\beta)\ln(\lambda)+\Ent(\alpha)+\Ent(\beta) +
        d\left[(1-\beta)\Ent({\alpha\over 1-\beta}) -
          \Ent(\alpha)\right]\right]\right\}\\
    & = & \exp (\Phi_1(\a,\b)\cdot n).
  \end{eqnarray*}
\end{proofof}

\begin{proofof}{Lemma~\ref{l:phi1}}
  We start by analyzing the
  first and second derivatives of $\Phi_1$. We refer to Claim~2.2
  and its proof in~\cite{DFJ} to establish the following:
  \begin{proposition}
    \label{p:derivs}
    \begin{romanenumerate}
    \item The function $\Phi_1$ has no local maxima on the boundary of $\cT$ and at
      least one local maximum in the interior of $\cT$.
    \item Any stationary point $(\a,\b)$ of $\Phi_1$ satisfies $\b=f(\a)$ and
      $\a=f(\b)$, where
      \begin{equation}
        \label{eq:rec_f}
        f(x) = (1-x)\left[1-\left({x\over \l(1-x)}\right)^{1/d}\right].
      \end{equation}

    \item All local maxima of $\Phi_1$ satisfy $\a+\b+d(d-2)\a\b \le 1$.
    \end{romanenumerate}
  \end{proposition}  
  
  We now claim that the function~$f$ defined in~(\ref{eq:rec_f}) also
  describes the relationship between the probabilities of occupancy on
  neighboring vertices in the infinite regular tree~$\Bethe^d$. Specifically,
  let $\mu$ be a simple semi-translation invariant Gibbs measure for the
  hard-core model with activity~$\lambda$ on~$\Bethe^d$. For any vertex~$v$ on
  the {\em even} (respectively {\em odd}) partition of the tree let $p_E$
  (respectively $p_O$) stand for $\mu[\mbox{$v$ is occupied}]$. (The definition
  does not depend on the choice of~$v$ since $\mu$ is semi-translation
  invariant.)  

  We will show below that $p_E$ and $p_O$ must satisfy the same
  relationship as in part~(ii) of Proposition~\ref{p:derivs}, i.e., it must be
  the case that $p_E=f(p_O)$ and $p_O=f(p_E)$. Indeed, recursive relationships
  similar to~$f$ for the probability of occupancy on the regular tree have been
  studied in the analysis of Gibbs measures on trees. See, e.g.,
  \cite{kelly,MSW04}. To see that $p_E=f(p_O)$ notice that by definition of the
  hard-core model, for any Gibbs measure~$\mu$ and every vertex~$v$,
  $\mu[\mbox{$v$ is occupied}] = {\l\over 1+\l}\cdot\mu[\mbox{all the neighbors
    of~$v$ are unoccupied}]$. On the other hand, if~$\mu$ is simple then
  $$\mu[\mbox{all the neighbors of~$v$ are unoccupied}] \;=$$
  $$\mu[\mbox{$v$ is occupied}] \;+\; \mu[\mbox{$v$ is unoccupied}]\prod_{i=1}^d
  \mu[\mbox{$u_i$ is unoccupied} \,|\, \mbox{$v$ is unoccupied}],$$
  where the
  $u_i$ are the neighbors of~$v$. Thus, if we let $\pcond_E= \mu[\mbox{$v$ is
    occupied} \,|\, \mbox{$u$ is unoccupied}]$, where $\mu$ is the
  semi-translation invariant measure under consideration and $v,u$ are two
  arbitrary vertices connected by an edge with $v$ even and $u$ odd then we get
  $p_O= {\l\over 1+\l} (p_O+(1-p_O)(1-\pcond_E)^d)$, i.e., 
$\pcond_E=1-  p_O ^{1/d} /
    \big(\l(1-p_O)\big)^{1/d}$. Plugging the latter expression for $\pcond_E$ into the
  trivial equation $p_E=(1-p_O)\pcond_E$ gives $p_E=f(p_O)$ as required. A
  similar derivation with the roles of even and odd vertices reversed gives
  $p_O=f(p_E)$.
  
  Since $f(x)$ is decreasing in~$x$, it follows that for all $\l$ 
  and $d$ there is
  a unique $p^*\in[0,1]$ such that $p^*=f(p^*)$, which means there is always a
  unique simple translation invariant Gibbs measure on the infinite tree. However,
  the analysis in~\cite{kelly} (see also~\cite{MSW04}) shows that for $\l\le
  \l_c(d)=  (d-1)^{d-1}/ (d-2)^d $ this measure is also the unique simple {\em
    semi}-translation invariant measure while for $\l>\l_c$ there are two more
  measures of the latter kind. This means that for $\l\le \l_c$, $\a=\b=p*$ is the
  unique solution to the system of equations $\a=f(\b)$ and $\b=f(\a)$ while for
  $\l>\l_c$ there exists $p_1<p^*<p_2$ such that $(\a=p_1,\b=p_2)$ and
  $(\a=p_2,\b=p_1)$ are (the only) two additional solutions. It is also easy to verify
  that $p^*=1/d$ at the critical activity~$\l_c$. Notice that $f$ is increasing
  in $\l$ and therefore so must be $p^*$, i.e., for $\l>\l_c$, $p^*>1/d$. By
  part~(iii) of Proposition~\ref{p:derivs} we get that for $\l>\l_c$ the point
  $\a=\b=p^*$, although stationary for $\Phi_1$, is not a local maximum. We
  conclude that for $\l>\l_c$ the maximum of $\Phi_1$ over the triangle~$\cT$ is
  achieved at $(\a=p_1,\b=p_2)$ and $(\a=p_2,\b=p_1)$.
\end{proofof}

\section{Logarithm of the second moment} \label{sec:logsec}
We would like to show concentration of $Z_G^{\alpha,\beta}$, at least for some
$\alpha,\beta$, by using the second moment. 
We begin with the function describing the exponent. We introduce 
the overlap   parameters $\gamma$ and $\delta$ as follows. We calculate the
contribution of pairs of independent sets such that both independent sets have
$\alpha n$ vertices on the right and $\beta n$ vertices on the left, and the
overlap on the right is $\gamma n$ and the overlap on the left is $\delta n$.
The constant in the exponent of this sum is then:
\begin{eqnarray*}
\label{eq:second_mom}
\Phi_2(\alpha,\beta,\gamma,\delta) &=& 2(\alpha+\beta)\ln(\lambda) \;+\\
& &
H(\alpha) + \alpha H({\gamma\over \alpha}) + (1-\alpha)H({\alpha-\gamma\over 1-\alpha}) \;+\;\\
& & H(\beta) + \beta H({\delta\over\beta}) + (1-\beta)H({\beta-\delta\over
  1-\beta}) \;+\;\\
& &
d\cdot\Psi_2(\alpha,\beta,\gamma,\delta,\eps),
\end{eqnarray*}
where $\Psi_2$ is the logarithm of the probability that two independent sets of the
structure   with parameters $(\alpha,\beta,\gamma,\delta)$   remain independent when a random matching is added. This is given by:
\begin{eqnarray*}
\label{eq:pr_two_sets}
\Psi_2(\a,\b,\g,\d,\e) & = & (1-2\b+\d)H({\g\over 1-2\b+\d}) - H(\g) \;+\\
& & (1-2\b+\d-\g)H({\epsilon\over 1-2\b+\d-\g}) +
      (\b-\d)H({\a-\g-\e\over \b-\d}) - (1-\g)H({\a-\g\over 1-\g}) \;+\\
& & (1-\b-\g-e)H({\a-\g\over 1-\b-\g-\e}) - (1-\a)H({\a-\g\over 1-\a}),
\end{eqnarray*}
where the three lines correspond to the probability of the following three
events respectively. The first is that the $\g n$ edges connected   to the
intersection
on the right avoid both sets on the left. The second is, conditioned on the first
occurring, that the $(\a-\g) n$ edges connected the first set on the right but
not
to the second avoid the first set on the left, where we sum over the    number $\e
n$ of edges that avoid both sets (and therefore take the maximum over $\e$). The
last event is, conditioned on the first two occurring, that the $(\a-\g)n$ edges
connected to the second set on the right but not to the first avoid the second
set on the left.

We would like to show that   $\Phi_2 \le 2\Phi_1$ given 
in Proposition~\ref{p:first_mom} for all relevant $\alpha,\beta$
and all $\gamma,\delta,\e$ (where the latter quantities have to make sense, e.g.,
$\g\le \a$, etc.), i.e., we want to show that:
\begin{eqnarray} \nonumber
\Gamma(\a,\b,\g,\d,\e) &:=& 
H(\a)+H(\b)-
\left[\a H({\g\over \a}) + (1-\a)H({\a-\g\over 1-\a}) + \b H({\d\over\b}) + (1-\b)H({\b-\d\over
  1-\b})\right]
   +\\ \label{eq:2nddiff}
 & & d \cdot [2\Psi_1(\a,\b)-\Psi_2(\a,\b,\g,\d,\e)] \ge 0\enspace .
\end{eqnarray}

\begin{lemma} \label{lem:gen_stat}
For all $\a,\b$ the point $(\a,\b,\g^{\ast},\d^{\ast},\e^{\ast})$
is a stationary point of both $\Phi_2$ and $\Psi_2$, 
where $\g^{\ast} = \a^2, \d^{\ast} = \b^2$
and $\e^{\ast} = \a(1-\a-\b)$.
\end{lemma}


\begin{lemma} \label{lem:gen_stat_zero}
  For all $\a,\b$ the point $(\a,\b,\g^{\ast},\d^{\ast},\e^{\ast})$ satisfies
  $\Gamma(\a,\b,\g^{\ast},\d^{\ast},\e^{\ast}) = 0$.
\end{lemma}

The meaning of the above two (easily verifiable) lemmas is that in order to show that
$\Gamma(\a,\b,\g,\d,\e)\ge 0$ for given $\a,\b$ and all $\g,\d,\e$ we have to
show that the stationary point $(\g^{\ast},\d^{\ast},\e^{\ast})$ is in fact the
global maximum of $\Phi_2$. While we believe this to be true for all relevant
$\a,\b$, we carried out the detailed calculations only for $\a,\b$ close to
the critical values.

\section{Calculation of the $2$nd moment around the critical point for $d \geq 3$}
\label{sec:ratio}
Here we show that the second moment is tight when
$d\geq3$ and $\alpha=\beta=1/d$. 

\subsection{The logarithm of second moment}
Recall that
\[
H(x) := -x \log (x)-(1-x) \log (1-x).
\]
Let:
\[
H_1(x,y) :=
  -x (\log (x)-\log (y)) + (x-y) (\log (y-x)-\log (y)).
\]
Then we have:
\begin{eqnarray*}
  \Psi_2(\a,\b,\g,\d,\e) = &\,&
  H_1(\gamma,(1-2*\beta+\delta)) - H(\gamma) +
  H_1(\e,(1-2*\beta+\delta-\gamma)) + \\ &\,&
  H_1((\alpha-\gamma-\e),(\beta-\delta)) -
  H_1((\alpha-\gamma),(1-\gamma)) + \\ &\,&
  H_1((\alpha-\gamma),(1-\beta-\gamma-\e))-
  H_1((\alpha-\gamma),(1-\alpha)):
\end{eqnarray*}
Therefore the important part of the log of the second moment is given by:
\begin{eqnarray*}
  f(\a,\b,\g,\d,\e) = &\,&
  2*(\alpha+\beta)*\ln(\lambda)     +
  H(\alpha)+H_1(\gamma,\alpha)+ \\
  &\,&
  H_1((\alpha-\gamma),(1-\alpha))
  +H(\beta)+H_1(\delta,\beta)+ H_1(\beta-\delta , 1-\beta)
  +d*\Psi_2(\a,\b,\g,\d,\e).
\end{eqnarray*}

We now find all stationary points in the region defined by all variables being
nonnegative as well as following constraints: \bel{region}
\alpha-\gamma-\epsilon \ge 0,\quad \beta-\delta\ge 0,\quad
1-2\beta+\delta-\gamma-\eps\ge 0.  \ee

\begin{lemma} \label{lem:exterior}
  Let $\alpha = \beta = 1/d$. Then the function $f$ as a function 
  of $(\g,\d,\e)$ obtains its maximum 
  the interior of the region~(\ref{region}). 
\end{lemma}

\begin{proof}
We know that for $\a,\b$, the independent point
$(\a,\b,\g^{\ast},\d^{\ast},\e^{\ast})$, where $\g^{\ast} = \a^2, \d^{\ast} =
\b^2$ and $\e^{\ast} = \a(1-\a-\b)$ is a stationary point of $f$ and we would
like to show that it is the global maximum for $\a = \b = 1/d$. 

Recall that 
$$
{\partial H_1(x,y)\over \partial x} = \log ({y-x\over x})
$$
and
$$
{\partial H_1(x,y)\over \partial y} = \log ({y\over y-x})
$$
so the derivatives of~$f$ are:
\[
\exp \left( \frac{df}{d\g} \right) = \,{\frac { \left(1-2\b+\d
      -\g-\e\right)^{d} \left(\alpha-\gamma-\epsilon \right) ^{d}
    \left(1-2\alpha+\gamma\right)^{d-1}}
  {\left(1-\beta-\gamma-\epsilon \right)^{d} \left(\beta-\delta-(\alpha-\gamma-\epsilon)\right)^{d} \left(\a-\g\right)^{d-2} \gamma}}
\]
\[
\exp \left( \frac{df}{d\d} \right) = \,{\frac {\left(\b-\d -
      (\alpha-\gamma-\epsilon)\right)^{d} \left(1-2\beta+\delta\right)^{d-1}}
  {\left(1-2\b+\d-\g-\epsilon\right)^{d} \left(\beta-\delta\right)^{d-2} \delta}}
\]
\[
\exp \left( \frac{df}{d\e} \right) = \,{\frac {
    \left(1-2\b+\d-\g-\epsilon\right)^{d} \left(\alpha-\gamma-\epsilon
    \right)^{d} \left(1-\alpha-\beta-\epsilon\right)^{d}}
  {\left(1-\beta-\gamma-\epsilon\right)^{d}
    \left(\b-\d- (\alpha-\gamma-\epsilon)\right)^{d}{\epsilon}^{d} }}   .
\]

For the proof of the lemma note that the derivatives (at least w.r.t. to one of the variables) go
to infinity (+ or -, in the right direction) as we approach any boundary point
of the region defined by~(\ref{region}).

The first derivative goes to $\infty$ as $\gamma \to 0$, the second as $\delta \to 0$ and the third as $\eps \to 0$. 
Similarly, the second derivative goes to $-\infty$  
as $\delta \to \beta$ the first as 
$\epsilon + \gamma - \delta \to 1 - 2 \beta$, and the third as 
$\gamma + \eps \to \alpha$.
This implies that the global maximum must be obtained at an interior point. 
\end{proof}

In order to proceed we need the second derivatives of $f$. 
\begin{lemma} \label{lem:f''}
\begin{eqnarray*}
\frac{\partial f}{\partial^2 \g} &=& -{d\over 1-2\b+\d-\g-\e} - {d\over
  \a-\g-\e} + {d-1\over 1-2\a+\g} -{d\over \b-\d-(\a-\g-\e)} + {d\over
  1-\b-\g-\e} +{d-2\over \a-\g} - {1\over\g}\\
\frac{\partial f}{\partial \d\partial \g} &=& {d\over 1-2\b+\d-\g-\e} +
{d\over \b-\d-(\a-\g-\e)}\\
\frac{\partial f}{\partial \e\partial\g} &=& -{d\over 1-2\b+\d-\g-\e} - {d\over
  \a-\g-\e} -{d\over \b-\d-(\a-\g-\e)} + {d\over
  1-\b-\g-\e}\\
\frac{\partial f}{\partial^2 \d} &=& -{d\over \b-\d-(\a-\g-\e)} + {d-1\over
  1-2\b+\d}-{d\over 1-2\b+\d-\g-\e}  + {d-2\over \b-\d} - {1\over \d}\\
\frac{\partial f}{\partial \e\partial \d} &=& {d\over \b-\d-(\a-\g-\e)} +
{d\over 1-2\b+\d-\g-\e}\\
\frac{\partial f}{\partial^2 \e} &=& -{d\over 1-2\b+\d-\g-\e} - {d\over
  \a-\g-\e} -{d\over 1-\a-\b-\e}
 + {d\over 1-\b-\g-\e} - {d\over \b-\d-(\a-\g-\e)} - {d\over\e}
\end{eqnarray*}
\end{lemma}

\begin{proof}
\end{proof}

Next we find for fixed $\d$ and $\g$ the $\e$ which maximizes $f$. 

\begin{lemma} \label{lem:esol}
For fixed values of $\a,\b,\d$ and $\g$ the maximum of $f$ is obtained for 
\[
\hat{\e}(\g,\d)\equiv \hat{\e}(\a,\b,\g,\d) = {1\over 2} \left[1+\a-\b-2\g -
\sqrt{(1-\a-\b)^2 + 4(\a-\g)(\b-\d)}\right]. 
\]
\end{lemma}

\begin{proof}
Solving for $\e$ in the equation ${\partial f\over \partial \e} = 0$ gives a
quadratic equation whose unique solution in the legal range of $\e$ is:
$\hat{\e}(\g,\d)\equiv \hat{\e}(\a,\b,\g,\d) = {1\over 2} \left[1+\a-\b-2\g -
\sqrt{(1-\a-\b)^2 + 4(\a-\g)(\b-\d)}\right]$. 
Note that $\hat{\e}$ is a maximizer since
${\partial f\over \partial^2 \e} < 0$ throughout the region~(\ref{region}). 
\end{proof}

We now define the function 
\[
g(\g,\d)\equiv g(\a,\b,\g,\d) = f(\g,\d,\hat{\e}(\g,\d)). 
\]

The proof will proceed by showing that the characteristic polynomial
of the Hessian matrix of~$g$ has only negative roots throughout the region
defined by 
\bel{region2}
0\le \g \le \alpha,\quad 0\le \delta \le \b.  
\ee

\begin{lemma} \lab{lem:g''}
The derivatives of~$g$ are:
\begin{eqnarray*}
{\partial g\over \partial\g}(\g,\d) & = & {\partial f\over \partial
  \g}(\g,\d,\hat{\e}(\g,\d)) \\
{\partial g\over \partial\d}(\g,\d) & = & {\partial f\over \partial
  \d}(\g,\d,\hat{\e}(\g,\d)) .
\end{eqnarray*}
The second derivatives are:
\begin{eqnarray*}
{\partial g\over \partial^2\g}(\cdot,\cdot) & = & {\partial f\over \partial^2
  \g}(\cdot,\cdot,\hat{\e}) + {\partial \hat{\e}\over \partial \g} \cdot
  {\partial f\over \partial\g\partial \e}(\cdot,\cdot,\hat{\e})\\
{\partial g\over \partial\d\partial\g}(\cdot,\cdot) & = & {\partial f\over \partial\d\partial
  \g}(\cdot,\cdot,\hat{\e}) + {\partial \hat{\e}\over \partial \g} \cdot
  {\partial f\over \partial\d\partial \e}(\cdot,\cdot,\hat{\e})\\
{\partial g\over \partial^2\d}(\cdot,\cdot) & = & {\partial f\over \partial^2
  \d}(\cdot,\cdot,\hat{\e}) + {\partial \hat{\e}\over \partial \d} \cdot
  {\partial f\over \partial\d\partial \e}(\cdot,\cdot,\hat{\e}).
\end{eqnarray*}
where
\begin{eqnarray*}
{\partial \hat{\e}\over \partial\g} &=& -1 + {\b-\d\over \sqrt{(1-\a-\b)^2+4(\a-\g)(\b-\d)}}\\
{\partial \hat{\e}\over \partial\d} &=& {\a-\g\over \sqrt{(1-\a-\b)^2+4(\a-\g)(\b-\d)}}.
\end{eqnarray*}
\end{lemma}

\begin{proof}
The calculation are straightforward, noting that ${\partial f\over \partial\e}(\cdot,\cdot,\hat{\e}) = 0$. 
\end{proof}

The major technical challenge is to prove that:
\begin{lemma} \label{lem:interior}
  Let $d \geq 3$ and $\alpha=\beta=1/d$. 
 Then the function $f$ has a unique stationary 
  point in
  the interior of the region~(\ref{region}). 
  This point is
\[
\begin{array}{ll}
\g^{\ast}=\d^{\ast}=\frac{1}{d^2},& 
\eps^{\ast}=\frac{1}{d} \left(1 - \frac{2}{d}\right) 
\end{array}
\]
and it is the maximum of the function.
\end{lemma}

Lemma~\ref{lem:interior} is proven in Lemma~\ref{lem:interior4} and
Lemma~\ref{lem:interior3} below. 

\begin{lemma} \label{lem:interior4}
  Let $d \geq 4$ and $\alpha=\beta=1/d$. 
 Then the function $f$ has a unique stationary 
  point in
  the interior of the region~(\ref{region}). 
  This point is. 
\[
\begin{array}{ll}
\g^{\ast}=\d^{\ast}=\frac{1}{d^2},& 
\eps^{\ast}=\frac{1}{d} \left(1 - \frac{2}{d}\right) 
\end{array}
\]
and it is the maximum of the function.
\end{lemma}

\begin{proof}
Given Lemma~\ref{lem:exterior} it suffices to show there is a unique local 
maximum of $f$ in the interior of~(\ref{region}). Given Lemma~\ref{lem:esol} 
it suffices to show that  
that~$g$ has a unique local maximum in the region~(\ref{region2}).
Using Lemma~\ref{lem:g''} for $\alpha=\beta=1/d$ we obtain:
\begin{eqnarray}
{\partial g\over \partial^2\d}  
& = & {\partial f\over \partial^2
  \d}(\cdot,\cdot,\hat{\e}) + {\partial \hat{\e}\over \partial \d} \cdot
  {\partial f\over \partial\d\partial \hat{\e}}(\cdot,\cdot,\hat{\e})\lab{dgddd} \\ 
& = & 
 -{d \over -\d+\g+\hat{\e}} + {d-1 \over
  1-2/d+\d}-{d\over 1-2/d+\d-\g-\hat{\e}}  + {d-2\over 1/d-\d} - {1\over \d}\non \\ &+& 
\left({d\over -\d+\g+\hat{\e}} + {d\over 1-2/d+\d-\g-\hat{\e}} \right)
{1/d-\g\over \sqrt{(1-2/d)^2+4(1/d-\g)(1/d-\d)}} \non \\
& = & 
 -{d \over -\d+\g+\hat{\e}} + {d-1 \over
  1-2/d+\d}-{d\over 1-2/d+\d-\g-\hat{\e}}  + {d-2\over 1/d-\d} - {1\over \d} \non\\ &+& 
\left({d\over -\d+\g+\hat{\e}} + {d\over 1-2/d+\d-\g-\hat{\e}} \right)
{1/d-\g\over 1 - 2 \g - 2 \hat{\e}}. \non
\end{eqnarray}

 
Note that for $d\ge 4$ we have 
\[
{1/d-\g\over 1 - 2 \g - 2 \hat{\e}} \leq 
{1/d-\g\over 1 - 2/d} \leq \frac{1}{d-2} \leq \frac{2}{d}.
\]
Therefore  
\begin{eqnarray*}
{\partial g\over \partial^2\d}  \le  
-{d-2\over \b-\d-(\a-\g-\hat{\e})} + {d-1\over
  1-2\b+\d}-{d-2\over 1-2\b+\d-\g-\hat{\e}}  + {d-2\over \b-\d} - {1\over \d}
\end{eqnarray*}
We now claim that the last expression is negative throughout the
region~(\ref{region2}). This follows since $\b-\d\ge
\b-\d-(\a-\g-\hat{\e})$, $1-2\b+\d \ge 1-2\b+\d-\g-\hat{\e}$ and $1-2\b+\d >
\d$. We therefore conclude that for $\a=\b=1/d$ the square derivative ${\partial g\over \partial^2\d} < 0$ throughout the region~(\ref{region2}). 
 
Since $g$ is symmetric in $\gamma$ and $\delta$ for $\alpha=\beta$ we conclude that ${\partial g \over \partial^2 \g} < 0$ throught the region as well.

We now note that ${\partial f\over \partial\g\partial\d} = {\partial
  f\over \partial\e\partial\d}$ and therefore we get that 
\bea
{\partial g\over \partial\d\partial\g} & = & 
\frac{\partial \hat\e}{\g} \cdot
 {\partial f\over \partial\e\partial
  \d}(\cdot,\cdot,\hat{\e})\lab{dgddg}\\
  &=&
{\b-\d\over \sqrt{(1-\a-\b)^2+4(\a-\g)(\b-\d)}} \cdot
 {\partial f\over \partial\e\partial
  \d}(\cdot,\cdot,\hat{\e}) \; > \; 0 \non
\eea
since ${\partial f\over \partial\e\partial
  \d} > 0$. Thus, since the square derivatives are negative, for any point
$(\g,\d)$, the 2nd derivative of $g$ along the line connecting the point
$(\g,\d)$ with its mirror $(\d,\g)$ is negative. (Notice that this 2nd
derivative is the sum of the square derivatives in each variable minus twice
the cross derivative.) Since $g$ is symmetric for $\a = \b = 1/d$, 
we conclude that the maximum of $g$ must be obtained on the line $\g=\d$. 

What remains to be shown is thus that the 2nd derivative of $g$ along the line
$\g=\d$ is negative. We proceed with $\a$ and $\b$ set equal to $1/d$.

This derivative equals to the sum of the square derivatives
plus twice the cross derivative. Since $g$ is symmetric and we are considering the
line $\g=\d$, it is enough to show that
$$
{\partial g\over \partial^2\d}(\d,\d) + {\partial g\over \partial\d\partial\g}(\d,\d) \; < \; 0$$
for every $0\le \d\le \b=1/d$. 
We use~\eqn{dgddd} and~\eqn{dgddg}, express the square root using the
formula in the statement of Lemma~\ref{lem:esol} 
 and note that in the present situation (where $\d=\g$)
$$
{\partial f\over \partial\e\partial
  \d} 
  =
  {\partial f\over \partial\e\partial
  \g},\quad  \mbox{ and also } \quad  
 1+ {\partial \hat{\e}\over \partial\g}  =   
{\partial \hat{\e}\over \partial\d}
= \frac{1/d-\d}{1-2\d-2 \hat \e}
$$ 
using Lemma~\ref{lem:g''}.
Thus, evaluating at $\e=\hat\e$ and $\g=\d$,
$$
{\partial g\over \partial^2\d}(\d,\d) + {\partial g\over \partial\d\partial\g}(\d,\d) =
{\partial^2 f\over \partial\d^2  }
+ 2
{\partial^2 f\over \partial\e\partial\d  } 
\cdot
{\partial \hat  \e\over \partial\d  }.
 $$
We find after some algebra that
\bea 
{\partial^2 f\over \partial\d^2  } &=&
\frac{-(d-2)P(d,\d,\hat \e)}{(d-2 -\hat \e d) \hat \e (d-2+d \d ) ( 1-d \d) \d},\non\\
 2
{\partial^2 f\over \partial\e\partial\d  } 
\cdot
{\partial \hat \e\over \partial\d  } &=&
\frac{2 (d-2) ( 1-d \d) }{ \hat \e(d-2 -\hat \e d)  ( 1- 2\d - 2\hat \e)}
\lab{second}
\eea 
where
\bean
P(d,\d,\hat \e) &=&  -\d d^3\hat \e-d^3 \d^2-d^3\d ^3+\d d^3 \hat \e^2+3 d^2 \d^2+2 d^2 \d \hat \e+\d d^2-\hat \e^2 d+\hat \e d-2 d \d-2 \hat \e\\
&=&  
 d ( \d d^2-1)\hat \e^2
-(d-2) ( \d d^2-1) \hat \e
+
d  \d ( 1-d  \d) (d-2+d  \d ).           
\eean
We may omit   the common factors $ (d-2)/( \hat \e(d-2 -\hat \e d))$. Note that, according to~\eqn{region}, $ \d + \hat \e \le \a=1/d$, and  hence we see that the expression in~\eqn{second} is positive and may be bounded above using $  (1-2/d)$ in place of $(1- 2 \d - 2\hat \e )$. Thus, it is enough to show
\bel{toshow}
P(d,\d,\hat \e)(1-2/d)
 -  2 \d (d-2+d \d ) ( 1-d \d)^2 >0.
\ee
Note that the formula for $\hat \e$ gives
$ 
\hat \e = \frac12 (a + b -\sqrt{a^2+b^2})
$ 
where $a=1-2/d$ and $b=2/d-2\d$. Realizing the numerator gives
$$ 
\hat \e = \frac{\frac{1}{2} (a + b +\sqrt{a^2+b^2}) (a + b -
  \sqrt{a^2+b^2})}{a+b+\sqrt{a^2+b^2}}
= \frac{ab} { a + b +\sqrt{a^2+b^2}}
= \frac{2(1-2/d)(1/d-\d)} { a + b +\sqrt{a^2+b^2}}
$$ 
Since both $a$ and $b$ are nonnegative, the denominator is at most $2(a+b) = 2(1-2\d)$, and so
$$ 
\hat \e \ge \hat \e_{{\rm min}}:= \frac{ (1-2/d)(1/d-\d)} { 1-2\d}.
$$
For an upper bound we will use, as  from~\eqn{region}, $\hat \e\le 1/d-\d$.

Now consider the parts of $P(d,\d,\hat \e)$  involving $\hat \e$:
$$
  ( 1- \d d^2 )\big((d-2)  \hat \e-  d \hat \e^2  \big).
$$
The second factor is monotonically increasing  in $\hat \eps$ and positive up to $1/d$ (as $d\ge 4$). So for $ \d  <1/d^2$, when the first factor is positive, $P$ is bounded below by substituting $\hat \e= \hat \e_{{\rm min}}$. Similarly, for  $ \d  \ge 1/d^2$, we obtain a lower bound by substituting  $\hat \e =  1/d-\d$.

Substituting $\hat \e =  1/d-\d$ into the left hand side of~\eqn{toshow} gives
$$
\frac{( 1-d   \d)  (2 d^4 \d^3+Q_1(d,\d) )}{  d^2}
$$
where
$$
Q_1(d,\d)= 2d^3(d-3)  \d^2 -d(d-1)(d-2)\d +  (d-2)(d-3).
$$
This quadratic in $\d$ has its minimum at $\d=(d-1)(d-2)/4d^2(d-3)$, where its value is 
$$
-\frac{(d-1)^2(d-2)^2}{8d (d-3)}+  (d-2)(d-3)
$$
which is easily seen to be positive for $d\ge 4$. Thus for $d\ge 4$, we have~\eqn{toshow} and are done in the case $ \d  \ge 1/d^2$.

Next consider $ \d  < 1/d^2$. Substituting $\hat \e =  \hat \e_{{\rm min}}$ into the left hand side of~\eqn{toshow} gives 
$$
\frac{ ( 1-d \d ) Q_2(d,\d)}{  d^4 ( 1-2\d)^2}$$
where
\bean
Q_2(d,\d) \!\!\!& =  & \!\!\!8d^6\d^5 +4d^5(d-8)\d^4+ d^4(-6d^2+8d+32)\d^3
+ d^3(10d^2-20d-8)\d^2\\
&&\!\!\!+d(d-2)(d^3-9d^2+8d-4)\d+(d-1)(d-2)^2.
\eean
Rewrite this as
$$
Q_2(d,\d) = 8d^6\d^5 +\big(4d^5(d-8)\d^4+d^4(8d+32)\d^3\big)
+  d^3(10d^2-20d-8-6d^3 \d -3d^2+6d)\d^2 
$$
$$
 +\big(  3d^4( d -2)\d^2 + d(d-2)(d^3-9d^2+8d-4)\d+(d-1)(d-2)^2\big).
$$

We argue that each of the four terms here is nonnegative (and the last is strictly positive) for $ \d  < 1/d^2$. The first is immediate. The second is clearly nonnegative for $d\ge 8$; otherwise the first part is minimized by $\d= 1/d^2$ which makes the whole term nonnegative. For the third term, the big factor is at least $ 7d^2-20d-8 >0$ as $d\ge 4$. The last and longest term is quadratic in $\d$. Its derivative with respect to $\d$ is, using the upper bound $1/d^2$ for $\d$, easily seen to be negative for all $d\le 7$. So for such $d$ we may substitute  $\d = 1/d^2$, and noting
\bel{useful}
 d^3-9d^2+8d-4 =  d(d-8)(d-1)-4
\ee
this whole term becomes
$$
3(d-2) + (d-1)(d-2)(d-4)(d+1) -4(d-2)/d >0.
$$
On the other hand, to cover the case $d\ge 8$ for the last term, ignoring the (clearly nonnegative) $\d^2$ term gives a linear function with positive constant term, so it is positive provided 
$$
d(d-2)(d^3-9d^2+8d-4)\frac{1}{d^2}+(d-1)(d-2)^2 >0.
$$
This is obvious for $d\ge 8$, using~\eqn{useful}.

This completes the proof that for $d\ge 4$ and $\a=\b=1/d$ the function $f$ has
a unique stationary point which is the global maximum of~$f$.
\end{proof}

\begin{lemma} \label{lem:interior3}
  Let $d = 3$ and $\alpha=\beta=1/3$. 
 Then the function $f$ has a unique stationary 
  point in
  the interior of the region~(\ref{region}). 
  This point is. 
\[
\g^{\ast}=\d^{\ast}=\eps^{\ast}=\frac{1}{9} 
\]
and it is the maximum of the function.
\end{lemma}

Some of the algebra the proof of the lemma below was performed using MAPLE. In
particular, MAPLE was used in order to symbolically factor polynomials and
calculate resultants.
\begin{definition}
  Let $x = (x_1,\ldots,x_k)$ be a vector of variables. Let $p(x,y) =
  \sum_{i=0}^{d_1} p_i(x) y^i$ and $q(x,y) = \sum_{i=0}^{d_2} q_i(x) y^i$ be
  two polynomials of degree $d_1$ and $d_2$ in the variable $y$. The {\em
    resultant} of $p$ and $q$ with respect to $y$, $R(p,q;y)$ is the
  polynomial of $x$ defined by the determinant of the {\em Sylvester matrix}
  of the two polynomials.
  
  \noindent
  The {\em Sylvester matrix} of $p$ and $q$ is the $(d_1+d_2) \times (d_1+d_2)$
  matrix defined by:
  \[
  \left(
    \begin{matrix}
      p_{d_1} & \ldots  & p_0     & 0       & \ldots  & 0      & 0 \\
      0       & p_{d_1} & \ldots  & p_0     & 0       & \ldots & 0   \\
      0       &  \ddots & \ddots  & \ddots  & \ddots  & 0      & 0   \\
      0       &  \cdots &  0      & p_{d_1} & \cdots  &  p_0   & 0 \\
      0       &  \cdots &  0      & 0       & p_{d_1} & \cdots &  p_0   \\
      q_{d_2} & \ldots  & q_0     & 0       & \ldots  & 0      & 0 \\
      0       & q_{d_2} & \ldots  & q_0     & 0       & \ldots & 0   \\
      0       &  \ddots & \ddots  & \ddots  & \ddots  & 0      & 0   \\
      0       &  \cdots &  0      & q_{d_2} & \cdots  &  q_0   & 0 \\
      0       &  \cdots &  0      & 0       & q_{d_2} & \cdots &  q_0   \\
    \end{matrix}\right)
  \]
\end{definition}

In particular, in the proof of the lemma, we will often use the
following well known result:
\begin{fact} \label{fact:res}
  If $(x,y)$ is a root
  of both $p$ and $q$ then $x$ is a root of their resultant $R(p,q;y)$.
\end{fact}

\begin{proof}
Given Lemma~\ref{lem:exterior} it suffices to show there is a unique local 
maximum of $f$ in the interior of~(\ref{region}). 

  We take the derivatives of~$f$ to obtain:
  \[
  \exp \left( \frac{df}{d\g} \right) = \,{\frac { \left( \epsilon-1+2\,\beta-\delta+\gamma \right) ^{3} \left( -\alpha+\gamma+\epsilon \right) ^{3}\\
      \mbox{} \left( 1-2\,\alpha+\gamma \right) ^{2}}{ \left(
        \beta-\delta-\alpha+\gamma+\epsilon \right) ^{3} \left( \gamma-\alpha
      \right) \left( -1+\beta+\gamma+\epsilon \right) ^{3}\gamma}}
  \]
  \[
  \exp \left( \frac{df}{d\d} \right) = \,{\frac { \left( 2\,\beta-\delta-1 \right) ^{2} \left( \alpha-\gamma-\epsilon-\beta+\delta \right) ^{3}}{ \left( \epsilon-1+2\,\beta-\delta+\gamma \right) ^{3}\\
      \mbox{} \left( \beta-\delta \right) \delta}}
  \]
  \[
  \exp \left( \frac{df}{d\e} \right) = \,{\frac { \left( \epsilon-1+2\,\beta-\delta+\gamma \right) ^{3} \left( \alpha-\gamma-\epsilon \right) ^{3}\\
      \mbox{} \left( \alpha-1+\beta+\epsilon \right) ^{3}}{{\epsilon}^{3} \left( \alpha-\gamma-\epsilon-\beta+\delta \right) ^{3} \left( -1+\beta+\gamma+\epsilon \right) ^{3}}}   .
  \]
  In particular, in order for $(1/3,1/3,\g,\d,\e)$ to be a stationary
  point, we must have     equality between the numerator and the
  denominator in the expressions above. In other words, if
  $(1/3,1/3,\g,\d,\e)$ is a stationary point then it is a zero of the
  following three polynomials  
  \begin{eqnarray*}
    G &:= &- \left( \epsilon-1+2\,\beta-\delta+\gamma \right) ^{3} \left( \alpha-\gamma-\epsilon \right) ^{3} 
    \mbox{} \left( 2\,\alpha-\gamma-1 \right) ^{2}-\gamma \left( \alpha-\gamma-\epsilon-\beta+\delta \right) ^{3} \left( \alpha-\gamma \right)  \left( -1+\beta+\gamma+\epsilon \right) ^{3}
    \\D &:=& \left( 2\,\beta-\delta-1 \right) ^{2} \left( \alpha-\gamma-\epsilon-\beta+\delta \right) ^{3}- \left( \epsilon-1+2\,\beta-\delta+\gamma \right) ^{3} \left( \beta-\delta \right) \delta  
    \\
    E &:= &-\delta\,\alpha+\delta\,\gamma+\gamma-2\,\beta \gamma-{\gamma}^{2}+\epsilon-\beta\epsilon-2\,\gamma\epsilon-{\epsilon}^{2}-\alpha+2\,\beta\alpha+\alpha \gamma+\alpha\epsilon,
  \end{eqnarray*}
  where, for the last equation, we used the fact that $x^3 = y^3$ implies $x=y$  for real  $x$ and $y$.  
  
  When we substitute $\alpha=\beta=1/3$ into these equations we get:
  \begin{eqnarray*}
    G &=&- \left( \epsilon-1/3-\delta+\gamma \right) ^{3} \left( 1/3-\gamma-\epsilon \right) ^{3} 
    \mbox{} \left( -1/3-\gamma \right) ^{2}-\gamma \left( -\gamma-\epsilon+\delta \right) ^{3} \left( 1/3-\gamma \right)  \left( -2/3+\gamma+\epsilon \right) ^{3}\\
    D &= & \left( -1/3-\delta \right) ^{2} \left( -\gamma-\epsilon+\delta \right) ^{3}- \left( \epsilon-1/3-\delta+\gamma \right) ^{3} \left( 1/3-\delta \right) \delta
    \\
    E &= &-1/3\,\delta+\delta\,\gamma+2/3\,\gamma-{\gamma}^{2}+\epsilon-2\,\gamma\epsilon-{\epsilon}^{2}-1/9   .
  \end{eqnarray*}
  In order to proceed, we eliminate variables using Fact~\ref{fact:res}
  and calculating resultants. We let $R(G,E) = R(G,E;\delta)$ and
  $R(D,E) = R(D,E;\delta)$. Then
  \[
  R(G,E) = \,{\frac {1}{59049}}\, \left( 81\,{\gamma}^{2}{\epsilon}^{2}+81\,\gamma{\epsilon}^{3}-27\,{\gamma}^{2}\epsilon-27\,\gamma{\epsilon}^{2}+9\,{\epsilon}^{3}+3\,{\gamma}^{2}+9\,\gamma\epsilon-\gamma \right)  \left( 3\,\gamma+3\,\epsilon-1 \right) ^{3} \left( 3\,\gamma+3\,\epsilon-2 \right) ^{3}
  \]
  and
  \begin{align*}
    R(D,E)\, = \,{\frac {1}{729}}\, \left( 3\,\gamma+3\,\epsilon-1
    \right)  \left(
      243\,{\gamma}^{4}{\epsilon}^{2}+972\,{\gamma}^{3}{\epsilon}^{3}+1458\,{\gamma}^{2}{\epsilon}^{4}+972\,\gamma{\epsilon}^{5}+243\,{\epsilon}^{6}-81\,{\gamma}^{4}\epsilon-648\,{\gamma}^{3}{\epsilon}^{2}-1566\,{\gamma}^{2}{\epsilon}^{3} \right. \\
    \mbox{}-1512\,\gamma{\epsilon}^{4}-513\,{\epsilon}^{5}+9\,{\gamma}^{4}+108\,{\gamma}^{3}\epsilon+486\,{\gamma}^{2}{\epsilon}^{2}+756\,\gamma{\epsilon}^{3}+369\,{\epsilon}^{4}-6\,{\gamma}^{3}-45\,{\gamma}^{2}\epsilon-132\,\gamma{\epsilon}^{2}\\\left.
      \mbox{}-105\,{\epsilon}^{3}+{\gamma}^{2}+6\,\gamma\epsilon+9\,{\epsilon}^{2}
    \right)   .
  \end{align*}
  
  Note $\e+\gamma=2/3$ is impossible since $\eps+\gamma\le \alpha$ in~\eqn{region}.
  
  We can eliminate  $\eps+\gamma = 1/3 $ as follows.
  Substituting $\eps=1/3-\gamma$ into equation $G$ gives
  \[
  -{\frac {1}{2187}}\,\gamma \left( 3\,\delta-1 \right) ^{3} \left( -1+3\,\gamma \right)
  \]
  and the only zeros of this are at $\delta$ or $\gamma=1/3$, which are
  on the boundary.
  
  It thus suffices to consider zeros of
  \[
  R(G,E) :=  \,81\,{\gamma}^{2}{\epsilon}^{2}+81\,\gamma{\epsilon}^{3}-27\,{\gamma}^{2}\epsilon-27\,\gamma{\epsilon}^{2}+9\,{\epsilon}^{3}+3\,{\gamma}^{2}+9\,\gamma\epsilon-\gamma
  \]
  and
  \begin{align*}
    R(D,E) := \,243\,{\gamma}^{4}{\epsilon}^{2}+972\,{\gamma}^{3}{\epsilon}^{3}+1458\,{\gamma}^{2}{\epsilon}^{4}+972\,\gamma{\epsilon}^{5}+243\,{\epsilon}^{6}-81\,{\gamma}^{4}\epsilon-648\,{\gamma}^{3}{\epsilon}^{2}\\
    \mbox{}-1566\,{\gamma}^{2}{\epsilon}^{3}-1512\,\gamma{\epsilon}^{4}-513\,{\epsilon}^{5}+9\,{\gamma}^{4}+108\,{\gamma}^{3}\epsilon+486\,{\gamma}^{2}{\epsilon}^{2}+756\,\gamma{\epsilon}^{3}+369\,{\epsilon}^{4}-6\,{\gamma}^{3}-45\,{\gamma}^{2}\epsilon\\
    \mbox{}-132\,\gamma{\epsilon}^{2}-105\,{\epsilon}^{3}+{\gamma}^{2}+6\,\gamma\epsilon+9\,{\epsilon}^{2}  .
  \end{align*}
  We let $R  = R(R(G,E;\d),R(D,E;\d);\g)$,   and find that
  \begin{align*}
    R = \,-243\, \left( 9\,\epsilon-1 \right)  \left( 27\,{\epsilon}^{2}-9\,\epsilon+1 \right)  \left( 81\,{\epsilon}^{4}-81\,{\epsilon}^{3}-27\,{\epsilon}^{2}+12\,\epsilon-1 \right) \\
    \mbox{} \left(
      1296\,{\epsilon}^{4}-1917\,{\epsilon}^{3}+840\,{\epsilon}^{2}-97\,\epsilon+6 \right)  \left( 3\,\epsilon-1 \right) ^{2} {\epsilon}^{4}    .
  \end{align*}
  
  The zero $\eps=1/9$ will be investigated below. The case $\eps=1/3$
  is on the boundary (it forces $\gamma=0$ by~(\ref{region})) so is not
  of interest at present. Similarly $\e=0$ is on the boundary.
  The other factors are
  
  \[
  27\,{\epsilon}^{2}-9\,\epsilon+1
  \]
  \[
  81\,{\epsilon}^{4}-81\,{\epsilon}^{3}-27\,{\epsilon}^{2}+12\,\epsilon-1
  \]
  \[
  1296\,{\epsilon}^{4}-1917\,{\epsilon}^{3}+840\,{\epsilon}^{2}-97\,\epsilon+6
  \]
  These polynomials have no roots in the required range as can be
  verified by elementary calculus. In the case of the    two quartics,  the
  plots in Figure~\ref{fig1} might help convince the reader of this.
  \figone
  
  So now we consider the only remaining case, $\eps=1/9$. Substituting
  $\e=1/9$ in $D$ and $G$ we obtain 
  \[
  D := \, \left( -1/3-\delta \right) ^{2} \left( -\gamma-1/9+\delta \right) ^{3}- \left( -2/9-\delta+\gamma \right) ^{3} \left( 1/3-\delta \right) \delta
  \]
  and
  \begin{equation}\label{Gdef}
    G := -(-2/9+\d+\g)^3(2/9-\g)^3(-1/3-\g)^2 - \g(-\g-1/9+\d)(1/3-g)(-5/9+\g)^3   .
  \end{equation}
  Letting $R=R(D,G,\g)$ we   find that $R$ factors as  
  $$ 
  {\frac {1}{984770902183611232881}}\, \left( 9\,\delta-1 \right) \\
  \mbox{} \left( 1458\,{\delta}^{3}+405\,{\delta}^{2}+24\,\delta+1
  \right)  
  \mbox{}  \left( 3\,\delta-1 \right) ^{3}P(\delta) 
  $$
  where  
  \begin{align*}
    P(\delta) = 1549681956\,{\delta}^{11}+2970223749\,{\delta}^{10}-157837977\,{\delta}^{9}-36669429\,{\delta}^{8} 
    \mbox{}+42830208\,{\delta}^{7}\\
    -35446896\,{\delta}^{6}-4331961\,{\delta}^{5}+1160487\,{\delta}^{4}+22734\,{\delta}^{3}+47529\,{\delta}^{2} 
    \mbox{}+12720\,\delta+64.
  \end{align*}
  
  We will consider the case $ \delta = 1/9$ later.
  Note that $\delta= 1/3$ is on boundary. The other factors are:
  \[
  1458\,{\delta}^{3}+405\,{\delta}^{2}+24\,\delta+1
  \]
  which is clearly positive for $\delta>0$, and
  $P(\delta)$ which is positive as shown by the plots for different intervals  in Figure
  2.
  \figtwo 
  More formally, this can be verified as follows. We first observe
  that the polynomial
  \[
  -157837977\,{\delta}^{8}-36669429\,{\delta}^{7}-35446896\,{\delta}^{5}-4331961\,{\delta}^{4}+12720,
  \]
  is decreasing and positive in the interval $0 \leq \delta \leq 0.18$.   So $ P(\delta)$ is positive for all such $\delta$.  Then
  expanding $P(\delta+0.18)$ we obtain that
  \[
  P(\delta+0.18) \geq 3000 + 9000 \delta - 3*10^5 \delta^2 -
  1.5*10^6 \delta^3
  \]
  for $\delta \geq 0$ (all other monomials have positive
  coefficients). It is easy to verify that the polynomial on the
  right hand side is positive for $\delta \in [0,0.09]$. Finally,
  looking at the polynomial $P(\delta+0.27)$ we see that all   of its
  coefficients are positive.
  
  The conclusion is that any interior stationary point must satisfy $\e =1/9$,
  $\delta =1/9$. Substituting these into~(\ref{Gdef}) gives
  \[
  {\frac {1}{177147}}\, \left( -1+3\,\gamma \right)  \left( 9\,\gamma-1 \right)  \left( 6561\,{\gamma}^{5}-5832\,{\gamma}^{4}+1377\,{\gamma}^{3}+18\,{\gamma}^{2}-36\,\gamma+8 \right)   .
  \]
  Zeros are at the boundary ($\gamma = 1/3$) or the crucial value $\gamma =
  1/9$. The large factor has no zeros in the relevant range.  This follows since
  both polynomials $6561\,{\gamma}^{2}-5832\,{\gamma}+1377$ and
  $18\,{\gamma}^{2}-36\,\gamma+8$ are positive for all $\g \in [0,1/4]$.
  Moreover for all $\gamma$ it holds that
  $6561\,{\gamma}^{2}-5832\,{\gamma}+1377 \geq 80$ and $80 \gamma^3 +
  18\,{\gamma}^{2}-36\,\gamma+8$ is positive in the interval $[0.25,1/3]$. (See
  Figure~\ref{fig3}.)  \figthree

  
  In order to conclude we have to show that the point $\g=\d=\e=1/9$ is in fact
  a local maximum. For this we calculate the Hessian matrix of the function $f$
  at that point to obtain:
%
  \[
  \left[ \begin {array}{ccc} -{\frac {243}{4}}&{\frac {81}{2}}&-{\frac
        {243}{4}}\\\noalign{\medskip}{\frac {81}{2}}&-{\frac {81}{2}}&{\frac
        {81}{2}}\\\noalign{\medskip}-{\frac {243}{4}}&{\frac {81}{2}}&-{\frac
        {405}{4}}\end {array} \right] .
  \]
  The characteristic polynomial of this Hessian matrix is
  \[
  {x}^{3}+{\frac {405}{2}}\,{x}^{2}+{\frac {45927}{8}}\,x+{\frac {531441}{16}}
  \]
  and all roots are less than zero, so we have a local maximum here.
  This concludes the proof of Lemma~\ref{lem:interior}.

\end{proof}

We would now like to conclude that for $\alpha$ and $\beta$ close
to $1/d$ it holds that there is a unique maximum at the stationary
point $(\a,\b,\g^{\ast},\d^{\ast},\e^{\ast})$. By concavity one
obtains
\begin{lemma} \label{lem:unique_max}
  There exists $\nu > 0$ such that if $|\alpha-1/d| < \nu$ and $|\beta - 1/d| <
  \nu$ then the function $g(\g,\d,\e)= f(\a,\b,\g,\d,\e)$ has a unique
  stationary point in the interior of~(\ref{region}), and this point is its
  global maximum.
\end{lemma}

\begin{proof}
  First note that since the function $f$ is continuous it follows that for
  sufficiently small $\nu$ the maximum of the function $g$ cannot be obtained on
  the boundary of the region. Therefore $g$ has at least one local maximum.
  
  From the continuity of the derivatives of $f$ it follows that for $\nu$
  sufficiently small, all stationary points of $g$ have to be $\eta$-close to
  the point $(1/d^2,1/d^2,1/d(1-2/d))$. 
  Moreover $\nu$ may be chosen such that $g$ is
  concave down on the $\eta$ neighborhood of $(1/d^2,1/d^2,1(1-2/d))$.  
  However, this
  implies that $g$ has a unique stationary point and it is a maximum. The proof
  follows.
\end{proof}

\subsection{The ratio of second to first moment}
\label{subsec:ratio}
So far we have only dealt with the logarithms of the first and second
moments. In order to apply the second moment method we need to consider   
the ratio between the moments more precisely. Using the quadratic
behavior of the function $f$ around the stationary point we obtain:
\begin{theorem} \label{thm:ratio}
  Let $d \geq3$. Then there exists $\eta > 0$ such that if $|\a-1/d| <
  \eta$, $|\b-1/d| < \eta$, the number $n$ is sufficiently large and $\a n,\b n$ are
  integers then
  \begin{eqnarray}  
    {\E_{\RG(n,d)}[(Z^{\a,\b}_G)^2]\over \E_{\RG(n,d)}^2[Z^{\a,\b}_G]}  &=&  (1+o(1)) \cdot
    {1\over 2\pi} \cdot {1\over \a(1-\a)\b(1-\b)}
    \int_{-\infty}^{\infty} \int_{-\infty}^{\infty} d \g d \d  \nonumber
\\    &\,\,\,&
    \left[{1\over \sqrt{2\pi}} \cdot {1\over 1-\a-\b} \cdot
      \sqrt{{(1-\a)(1-\b)\over \a\b}} \int_{-\infty}^{\infty} d\e 
      \exp \left\{ \frac{1}{2d} (\g,\d,\e)  H_f (\g,\d,\e)^t \right\} \right]^d
 \label{eq:ratio}
  \end{eqnarray}
  where $H_f$ is the Hessian of the function $f$ at the point
  $\a,\b,\g^{\ast}=\a^2,\d^{\ast}=\b^2,\e^{\ast}=\a(1-\a-\b)$.
\end{theorem}

\begin{proofof}{Theorem~\ref{thm:ratio}}
  We use the approximation
  $${n \choose a n} = (1+o(1)) {1\over \sqrt{2\p n}}\cdot {1\over \sqrt
    {a(1-a)}}\cdot e^{H(a) n}, $$
  and thus
  $${bn \choose a n} = (1+o(1)) {1 \over \sqrt{2\p n}} \cdot{\sqrt{b \over
      a(b-a)}}\cdot e^{bH(a/b)n}.$$
  We have:
  \begin{eqnarray*}
    \E_\RG[(Z^{\a,\b}_G)] & = &
    \l^{\a+\b}{n\choose \a n}{n\choose \b n}
    \left[{{(1-\b)n\choose \a n}\over {n\choose \a n}}\right]^d
  \end{eqnarray*}
  and
  \begin{eqnarray*}
    \E_\RG[(Z^{\a,\b}_G)^2] & = & 
    \l^{2(\a+\b)}{n\choose \a n}{n\choose \b n} \sum_{\g,\d}
    {\a n\choose \g n}{(1-\a)n\choose (\a-\g)n}{\b n\choose \d n}{(1-\b)n\choose
      (\b-\d)n}\\
    &&
    \left[{{(1-2\b+\d)n\choose \g n}\over {n\choose \g n}}
      \sum_{\e} {{(1-2\b+\d-\g) n \choose \e n}
        {(\b-\d)n\choose (\a-\g-\e)n} \over {(1-\g) n\choose (\a-\g)n}}
      {{(1-\b-\g-\e)n \choose (\a-\g)n}\over {(1-\a)n\choose (\a-\g)n}}\right]^d
    \enspace.
  \end{eqnarray*}
  Thus,
  \begin{eqnarray*}
    {\E_\RG[(Z^{\a,\b}_G)^2]\over \E_\RG[Z^{\a,\b}_G]^2} & = &
    {n\choose \a n}^{-1}{n\choose \b n}^{-1} \left[{{(1-\b)n\choose \a n}\over
        {n\choose \a n}}\right]^{-2d} \\
    &&\sum_{\g,\d}
    {\a n\choose \g n}{(1-\a)n\choose (\a-\g)n}{\b n\choose \d n}{(1-\b)n\choose
      (\b-\d)n} \\
    &&
    \left[{{(1-2\b+\d)n\choose \g n}\over {n\choose \g n}}
      \sum_{\e} {{(1-2\b+\d-\g) n \choose \e n}
        {(\b-\d)n\choose (\a-\g-\e)n} \over {(1-\g) n\choose (\a-\g)n}}
      {{(1-\b-\g-\e)n \choose (\a-\g)n}\over {(1-\a)n\choose (\a-\g)n}}\right]^d
  \end{eqnarray*}
  Using the above approximations, we get
  \begin{eqnarray*}
    {\E_\RG[(Z^{\a,\b}_G)^2]\over \E_\RG[Z^{\a,\b}_G]^2} 
    & = & (1+o(1))\left({1\over \sqrt{2\p n}}\right)^{d+2}
    \a(1-\a)\b(1-\b) \left({1-\a-\b\over (1-\a)(1-\b)}\right)^d\\
    &&\sum_{\g,\d}
    {1\over (\a-\g)(\b-\d)}\sqrt{{1 \over \g\d(1-2\a+\g)(1-2\b+\d)}}\\
    && \;\;\;\;\;\;e^{[-H(\a) - H(\b) +\a H({\g\over\a}) + (1-\a)H({\a-\g\over 1-\a}) +
      \b H({\d\over \b}) + (1-\b)H({\b-\d\over 1-\b})]n}
    \\
    &&\;\;\;\;\;\;\left[\sqrt{(1-2\b+\d)(1-2\a+\g)(\b-\d)(\a-\g)}\right.\\
    &&\;\;\;\;\;\;\sum_{\e}
    \sqrt{{(1-\b-\g-\e) \over
        \e(1-2\b+\d-\g-\e)(\a-\g-\e)(\b-\d-\a+\g+\e)(1-\b-\a-\e)
      }}\\
    && \;\;\;\;\;\;\;\;\;\;\;\;
    \left. e^{[\Psi_2(\a,\b,\g,\d,\e)-2\Psi_1(\a,\b)]n}\right]^d \\
    & = & (1+o(1))\left({1\over \sqrt{2\p n}}\right)^{d+2}
    \a(1-\a)\b(1-\b) \left({1-\a-\b\over (1-\a)(1-\b)}\right)^d\\
    &&\sum_{\g,\d}
    {1\over (\a-\g)(\b-\d)}\sqrt{{1 \over \g\d(1-2\a+\g)(1-2\b+\d)}}\\
    &&\;\;\;\;\;\;\left[\sqrt{(1-2\b+\d)(1-2\a+\g)(\b-\d)(\a-\g)}\right.\\
    &&\;\;\;\;\;\;\sum_{\e}
    \sqrt{{(1-\b-\g-\e) \over
        \e(1-2\b+\d-\g-\e)(\a-\g-\e)(\b-\d-\a+\g+\e)(1-\b-\a-\e)
      }}\\
    && \;\;\;\;\;\;\;\;\;\;\;\;
    \left. e^{\frac{\Gamma(\a,\b,\g,\d,\e)}{d}n}\right]^d \\
  \end{eqnarray*}
  For $d \geq 3$ and $\a,\b$ close to $1/d$, it follows from
  Lemma~\ref{lem:unique_max} that the function $\Gamma(\a,\b,\g,\d,\e)$ has a
  unique maximum and that the function decays quadratically around this point.
  This implies that all terms in the sum above that have $\g,\d$ or $\e$ more
  than $\Omega(1)$ away from the maximal value have an exponentially low
  contribution.  Thus, up to losing a factor of $1+o(1)$, we can plug in the
  above values of $\g,\d,\e$ except into the exponential terms. This gives:
  \[
  (1+o(1)) \cdot
  {1\over 2\pi n} \cdot {1\over \a(1-\a)\b(1-\b)} 
  \sum_{\g,\d} \left[{1\over \sqrt{2\pi n}} \cdot {1\over 1-\a-\b} \cdot
    \sqrt{{(1-\a)(1-\b)\over \a\b}} \sum_{\e}
    e^{\frac{\Gamma(\a,\b,\g,\d,\e)}{d}n}\right]^d.
  \]
  Finally we use the quadratic approximation of $\Gamma(\a,\b,\g,\d,\e)$ around
  $\a,\b,\g^{\ast}(\a,\b),\d^{\ast}(\a,\b),\e^{\ast}(\a,\b)$ and the standard
  approximation of integral by sums to arrive at an integral formula.  For this
  we recall that
  \[
  \Gamma(\a,\b,\g^{\ast},\d^{\ast},\e^{\ast}) = 0
  \]
  and we note that the
  Hessian of $\Gamma$ equals the Hessian of $f$ as 
  $\Gamma(\a,\b,\g,\d,\e) = f(\a,\b,\g,\d,\e) + 2 \Phi_1(\a,\b)$ and 
  $\Phi_1$ does not depend on $\g,\d,\e$. We thus obtain~(\ref{eq:ratio}).
\end{proofof}

We are in a
position to give the proof of Theorem~\ref{thm:2nd}.

\begin{proofof}{Theorem~\ref{thm:2nd}}
In order to prove the theorem, 
we need to calculate the expression in~(\ref{eq:ratio}).   
  Using the derivatives calculated in Lemma~\ref{lem:f''} we have:
$$
  H_f =   \left( \begin{matrix} 
  h_{11} & h_{12}  &  h_{13}\\
  h_{12} & h_{22} & h_{23}      \\
  h_{13} &h_{23}  & h_{33}
  \end{matrix} \right)
 $$
where 
\bean
h_{11}&=&
\frac{\a+d-2}{\a(\a-1)^2}-\frac{\b+d\a}{ \a^2\b}+\frac{d}{(1-\a)(1-\b)}-\frac{d}{\b(1-\b)(1-\a-\b)},\\
h_{12}&=&  \frac{d}{\b(1-\b)(1-\a-\b)} \\
h_{13}&=& - \frac{d (1-\a-2\b+2\a\b+\b^2)}{\a\b(1-\a)(1-\b)(1-\a-\b)} \\
h_{22}&=&  - \frac{1-\a-\b+d\a\b}{\b^2( 1-\b)^2(1-\a-\b)}\\
h_{23}&=&   \frac{d}{\b(1-\b)(1-\a-\b)}\\
h_{33}&=&  - \frac{d(1-\a-\b+2\a\b)}{\a\b(1-\a)(1-\b)(1-\a-\b)}.
 \eean
  We let
 $$
  g(\g,\d,\e) = \frac{1}{2d} (\g,\d,\e) H_f(\g,\d,\e)^t
$$
and observe that this is quadratic in $\e$. Then putting
$$
  h(\g,\d) = 
 \left[{1\over \sqrt{2\pi}} \cdot {1\over 1-\a-\b} \cdot
      \sqrt{{(1-\a)(1-\b)\over \a\b}} \int_{-\infty}^{\infty} d \e \,
          e^{g(\g,\d,\e)} \right]^d 
$$
as required,  Gaussian integration gives  
  \[
  h(\g,\d) =
 \big(A_d  \exp \left( -B_d   \right)\big)^d,  
  \]
  where
 $$
A_d= \frac{(1-\a)(1-\b)}{\sqrt{(1-\a-\b+2\a\b) (1-\a-\b)}}
 $$
 and $B_d$ is much more complicated and is a quadratic polynomial in $\g$ and in $\d$.
In fact
$$B_d = \frac{B^2-4AC}{4A}$$
 where
\bean
  A  &=&  \frac{ 2\beta \alpha -\beta +1-\alpha }{2\beta \alpha(1-\beta -\alpha +\beta \alpha )(\beta +\alpha -1)},\\
B   &=& \frac{  -\alpha \gamma +2\beta \gamma \alpha +\gamma -2\beta \gamma +\beta ^2\gamma +\delta\alpha ^2-\delta\alpha }{\beta \alpha (-\alpha +\beta \alpha +1-2\beta +\beta ^2)(\alpha -1)},\\
 C &=& \frac{ \gamma ^2 (2d\beta \alpha ^3-d\alpha ^3+2d\beta ^2\alpha ^2-5d\alpha ^2\beta +2d\alpha ^2+d\beta ^3\alpha -3d\beta ^2\alpha -\beta ^2\alpha +3d\beta \alpha +\beta \alpha -d\alpha -\beta ^3+2\beta ^2-\beta )}{2d\alpha ^2(\beta -1)(-1+\beta +\alpha )(\alpha -1)^2\beta }
 \\
 &&+\frac{\delta \gamma }{2\beta (-\alpha +\beta \alpha +1-2\beta +\beta ^2)}  + 
\frac{(\gamma d\beta ^2-\gamma d\beta +d\beta \alpha \delta-\delta\alpha +\delta-\delta\beta )\delta}{2d\beta ^2(-\alpha +\beta \alpha +1-2\beta +\beta ^2)(\beta -1)}.
\eean

Integrating again, we obtain   
$$
{1\over 2\pi} \cdot {1\over \a(1-\a)\b(1-\b)}
\int_{-\infty}^{\infty} \int_{-\infty}^{\infty} d \g \, d \d \,
  h(\g,\d) =
  \tau^{\a,\b}(d)
$$
as required.
\end{proofof}

\section{Asymptotically almost sure results} \label{sec:aas}

In this section we prove Theorem~\ref{thm:aas} using the small graph
conditioning method. 

The small subgraph conditioning method has some chance of applying to
a random variable $Y$ when
the variance $\Var(Y)$ is of the same order as $(\E\,Y)^2$. This is
indeed the case for $Z_G^{\alpha,\beta}$ for $\alpha$ and $\beta$ such that
the conclusion of Theorem~\ref{thm:2nd} holds. In random regular graphs the only interesting local structures that occur with nonvanishing probability are short cycles. The method usually ``explains" the variance of $Y$ by the interaction between the numbers of short cycles and the random variable $Y$.  
For details, see \cite[Theorem 9.12--Remark 9.18]{JLR} and
\cite[Theorem 4.1, Corollary 4.2]{models}. 
The following is  a simplification of the latter.
 Here $[x]_m =x(x-1)\cdots(x-m+1)$  denotes the falling factorial and a.a.s.\ denotes ``asymptotically almost surely," i.e.\ the probability tends to 1 as $n\to\infty$.

\begin{theorem}\label{sscond}
Let $\lambda_i>0$ and $\delta> -1$ be real numbers for
$i=1,2,\ldots$. Let $\omega(n) \to 0$ 
and suppose that for each $n$ there are random variables $X_i =
X_i(n)$, $i=1,2,\ldots$
and $Y=Y(n)$, all defined on the same probability space $\cG = \cG_n$
such that $X_i$ is nonnegative integer valued, $Y$ is nonnegative and
$\E\,Y>0$ (for $n$ sufficiently large).  Suppose furthermore that
\begin{description}
\item{(i)} For each\, $k\geq 1$, the variables $X_1,\ldots ,X_k$ are asymptotically
independent Poisson random variables with $\E\,X_i\to\lambda_i$,
\item{(ii)}  for every finite sequence $m_1,\ldots ,m_k$ of nonnegative
integers,
\bel{delta}
\frac{\E(Y[X_1]_{m_1}\cdots [X_k]_{m_k})}{\E\,Y}
  \rightarrow \prod_{i=1}^k \big(\lambda_i(1+\delta_i)\big)^{m_i}
\ee

\item{(iii)} $\sum_i \,\lambda_i\,{\delta_i}^2 < \infty$,
\item{(iv)} $\E\,Y^2/(\E\,Y)^2\leq \exp(\sum_i \lambda_i\,{\delta_i}^2) + o(1)$
as $n\rightarrow\infty$.
\end{description}
Then $Y> \omega(n) \E Y$ a.a.s.
\end{theorem}

The probability space $\Omega_n$ we are working with is the set of bipartite (multi)graphs obtained by taking $d$ random perfect matchings between two sets $V_1$ and $V_2$ of $n$ vertices each. This probability space is contiguous with a uniformly random $d$-regular graph (see the note after  the proof of~\cite[Theorem 4]{MRRW}), and hence, once we have verified the hypotheses of the theorem, $Y>0$ is a.a.s.\   in the uniform model as well (as well as various other models contiguous to it).

Let 
\begin{equation} \label{eq:defY}
Y = \l^{-(\alpha+\beta) n}  Z_G^{\alpha,\beta} 
\end{equation}
be the number of independent sets with
 $\alpha n$ vertices in $V_1$ and $\beta n$ in $V_2$. 
Let the variable
 $X_i$ be the number of cycles in the graph of length $i$. 
(We will
 apply the theorem only for even integers $i$, which is valid by a
 trivial change of notation.)

\begin{theorem} \label{thm:YOK}
Let $\a$ and $\b$ be such that the conclusion of
Theorem~\ref{thm:2nd} holds.   
Then the random variables $Y$ defined in~(\ref{eq:defY}) and 
$X_i$ defined by the number of cycles of length $i$ (for even $i$), satisfy the conditions of Theorem~\ref{sscond} 
\end{theorem}

\begin{proofof}{Theorem~\ref{thm:aas}}
By Theorem~\ref{thm:YOK} and Theorem~\ref{sscond} it follows that
a.a.s. we have $Y \geq \frac{1}{n} \E Y$. This clearly implies that 
\[
Z_{G}^{\alpha,\beta} \geq \frac{1}{n} \E Z_G^{\alpha,\beta}
\] 
as needed.
\end{proofof}

  An alternative valid conclusion of the theorem (not as quoted above) in this application, is that, if we sample from the random graphs with weight proportional to the number of independent sets, then the model of random graphs we get is contiguous to the original: events that are a.a.s.\ true in one model are also a.a.s.\ true in the other.

We now prove {Theorem~\ref{thm:YOK}}

\begin{lemma} \label{lem:lam_i}
Condition (i) holds with 
\bel{lambda}
\lambda_i = \frac{ r(d,i)}{i}
\ee
where $r(d,i)$ is the number of ways one can properly edge color a
cycle of length $i$ with $d$ colors. 
\end{lemma}

\begin{proof}
This follows using the standard techniques (e.g.\ see Bollob{\'a}s
book or~\cite{models}). The reasoning goes as follows. 
There are asymptotically $n^i/i$
positions for the cycle to be in, and, given the perfect matchings
that the edges of a given cycle belong to ($r(d,i)$ possibilities),
the probability it occurs is easily seen to be asymptotic to $n^{-i}$.
\end{proof}

For part (ii), as usual we do a calculation for $\E(Y X_i ) /\E\,Y $, which determines $\delta_i$, and observe that the same calculation is easily extended to the arbitrary moments required for verifying this  part of the theorem with the value of $\delta$ so obtained.

\begin{lemma} \label{lem:del_i}
\[
\frac{\E(Y X_i )}{\E\,Y} \to \l_i(1+\delta_i),
\]
where
\[
\delta_i =\frac{ \a^{i/2}\b ^{i/2}}{(1-\a)^{i/2}(1-\b)^{i/2}}.
\]
\end{lemma}

\begin{proof}
Note that 
$$
\E(Y X_i ) =
 \sum_{S,T} \sum_C \pr(A_{S,T}\wedge A_{C})
$$
where $S$ and $T$ denote subsets of $V_1$ and $V_2$  of sizes $\alpha n$ and $\beta n$ respectively, $C$ denotes a possible position of a cycle (not joining any vertices of $S$ and $T$), $A_{S,T}$  is the   event that the random graph has $S\cup T$ as an independent set, and $A_{C}$ the  event that it contains $C$.
Similarly,
$$
\E(Y ) =
 \sum_{S,T}   \pr(A_{S,T}).
$$

We proceed with $i$ even. We have 
$$
\E(Y X_i ) =
\frac{1}{i}\sum_{S,T} \sum_\xi \sum_\eta P_1 P_2
$$
where 
\begin{description}
\item{$\bullet$}    The leading factor $1/i$ accounts for the fact that we will count cycles rooted at a vertex in $V_1$ (which can be done in $i/2$ ways) and oriented (2 ways),
\item{$\bullet$} $S$ and $T$ denote subsets of $V_1$ and $V_2$  of sizes $\alpha n$ and $\beta n$ respectively, 
\item{$\bullet$} $\xi$ denotes a proper $d$-edge-coloured  rooted, oriented  $i$-cycle ($r(d,i)$ possibilities), in which the vertices are 2-coloured, black and white, with no two black vertices adjacent. The color of the edges will prescribe which of the $d$ perfect matchings an edge of a (potential) cycle will belong to. The black vertices will prescribe which of the cycle vertices are   members of $S\cup T$.
\item{$\bullet$} $\eta$ denotes a position that an $i$-cycle  can be in (i.e.\ the exact vertices it traverses, in order) such that prescription of the vertex colors of $\xi$ is satisfied,
\item{$\bullet$} $P_1$ is the probability that a random graph in $\Omega$ contains a cycle $C$ in the given position $\eta$ with the edge colors prescribed by $\xi$  in accordance with which matchings contain the edges of $C$,
\item{$\bullet$} $P_2$ is the conditional  probability that the random graph respects $S\cup T$ as an independent set, given that it contains $C$ as in the definition of $P_1$.
\end{description}

Since all quantities concerned are independent of $S$ and $T$ (provided they have the correct cardinalities), we can fix $S$ and $T$ and write 
\bel{ratio}
\frac{\E(Y X_i )}{\E\,Y}=
\frac{1}{i} \sum_\xi \sum_\eta \frac{P_1 P_2}{P_3}
\ee
where $P_3$ is the probability that the random graph respects $S\cup
T$ as an independent set.

As noted before, $P_3$ is the probability that there are no
edges between $T$ and $S$ in each of the $d$ matchings. In other
words, 
\[
P_3 = 
\left( \frac{[(1-\beta)n]_{\alpha n}} {[n]_{\alpha n}} \right)^d.
\]
For $\ell =1$ and $2$ let $j_ \ell(\xi)$ denote the number of black
vertices in the coloring prescribed by $\xi$ that lie in $V_i$. 
Next we show that $P_2$ asymptotically depends only on
$j_{\ell}$. For $k=1,\ldots,d$, let $e(k)$ denote the number of edges
of color $k$ in $\xi$. Let $f_1(k)$ denote the number of edges of
color $k$ adjacent to black vertices of $S$ and $f_2(k)$ denote the
number of edges of color $k$ adjacent to black vertices in $T$. 
Then given $\xi$ that contains no edges
connecting $S$ and $T$, the probability that $S,T$ is an independent
set is given by:
\[
P_2 = \prod_{k=1}^d \frac{ [(1-\b) n - e(k) + f_2(k)]_{(\a n - f_1(k))}}
              { [n - e(k)]_{(\a n - f_1(k))} }
\]
In order to calculate the asymptotics of $P_2/P_3$ we observe that
\[
\frac{[n]_{\a n}}{[n-e(k)]_{(\a n - f_1(k))}} \sim 
((1-\a)n)^{f_1(k)} \frac{[n]_{\a n}}{[n-e(k)]_{\a n}} \sim
((1-\a)n)^{f_1(k)} (\frac{1}{1-\a})^{e(k)}
\]
and similarly
\[
\frac{[(1-\b)n]_{\a n}}{[(1-\b)n-e(k)+f_2(k)]_{(\a n - f_1(k))}} \sim 
((1-\a-\b)n)^{f_1(k)} (\frac{1-\b}{1-\a-\b})^{e(k)-f_2(i)}
\]
Therefore
\[
\frac{P_2}{P_3} \sim ((1-\a)n)^{\sum_{k} f_1(k)}
(\frac{1}{1-\a})^{\sum_i e(k)}
((1-\a-\b)n)^{-\sum_i f_1(k)} (\frac{1-\b}{1-\a-\b})^{\sum_k f_2(k)-e(k)}
\]
Recalling that $\sum_k f_1(k) = 2 j_1, \sum_k f_2(k) = 2 j_2$ and 
$\sum_k e(k) = i$ we obtain that 
\[
\frac{P_2}{P_3} \sim 
\frac{(1-\a-\b)^{i-2j_1-2j_2}}{(1-\a)^{i-2j_1}(1-\b)^{i-2j_2}}.
\]

Clearly 
$$
P_1\sim n^{-i} 
$$
and the number of terms in the summation over $\eta$ for which $S,T$
may be an independent set (i.e.\ number of possible $\eta$) is asymptotic to
$$
\alpha^{j_1}(1-\alpha)^{i/2-j_1}\b^{j_2}(1-\b)^{i/2-j_2}n^i.
$$
Thus~\eqn{ratio} is asymptotic to
\begin{equation} \label{eq:gfunc}
\frac{1}{i} \sum_{ j_1,j_2} a_{ i,j_1,j_2}x^iy^{j_1}z^{j_2}
\end{equation}
where $ a_{i,j_1,j_2}$ is the number of possible $\xi$ of length $i$ with $j_\ell$ black vertices in $V_\ell$ ($\ell =1$ and 2), and  
$$
x=   \frac{1-\a-\b}{\sqrt{(1-\a)(1-\b)}} ,\qquad
y = \frac{\a(1- \a)}{(1-\a-\b)^2}, \qquad
z= \frac{\b(1- \b)}{(1-\a-\b)^2}.
$$

Define the matrix
$$A=
\left[
\begin{array}{cccc}
0 & 0 & 0 & 1 \\
0 & 0 & y & 1\\
0 & 1 & 0 & 0 \\
z & 1 & 0 & 0
\end{array}
\right].
$$
Each entry of $A$ refers to a transition from one state to the next as
we traverse the cycle $C$. 
The first row and column refer to a black vertex in $V_1$, the second
to a white vertex in $V_1$, the third to a black vertex in $V_2$, and the fourth to a white vertex in $V_2$. Then the trace of $A^i$ counts the possible $\xi$ weighted by $y^{j_1}z^{j_2}$, except for the edge coloring, of which there are $r(d,i)$ possibilities. Hence~\eqn{ratio} is asymptotic to $r(d,i)x^i{\rm tr}(A^i)/i$, with $x$, $y$ and $z$ defined as above.
Letting $\mu_k$, $k=1\ldots 4$ denote the eigenvalues, a little computation gives
\bean
\mu_1^2 &=& \mu_2^2 = \frac12 \Big(u+\sqrt{u^2-4v}\Big),\\
\mu_3^2 &=& \mu_4^2 = \frac12\Big(u-\sqrt{u^2-4v}\Big),
\eean
where $u=1+y+z$ and $v=yz$. Substituting the values of $y$ and $z$ gives  
$$
\mu_1^2 = \frac{\a\b}{(1-\a-\b)^2},\qquad  \mu_3^2 = \frac{(1-\a)(1-\b)}{(1-\a-\b)^2}.
$$
Thus, recalling that $i$ is even,
 and ${\rm tr}(A^i) = 2\mu_1^i+2\mu_3^i$,
 we have from~\eqn{ratio}
\bean
\frac{\E(Y X_i )}{\E\,Y} & \sim & \frac{r(d,i)x^i}{i} \cdot \frac{ \a^{i/2}\b ^{i/2}}{(1-\a-\b)^i}
+ \frac{r(d,i)x^i}{i} \cdot \frac{(1-\a)^{i/2}(1-\b)^{i/2}}{(1-\a-\b)^i}\\
& = & \frac{r(d,i) }{i} \cdot \frac{ \a^{i/2}\b ^{i/2}}{(1-\a)^{i/2}(1-\b)^{i/2}}
+ \frac{r(d,i)}{i}.
\eean
Now~\eqn{lambda} gives that in~\eqn{delta},
$$
\delta_i =\frac{ \a^{i/2}\b ^{i/2}}{(1-\a)^{i/2}(1-\b)^{i/2}}.
$$
for even $i\ge 2$.
\end{proof}

Verification of (ii) for arbitrary sequences $m_1, \ldots $ is based
on a straightforward extension of the above argument which we sketch
briefly. 
\begin{lemma} \label{lem:deltas}
For every finite sequence $m_1,\ldots ,m_k$ of nonnegative
\[
\frac{\E(Y[X_2]_{m_1}\cdots [X_{2k}]_{m_k})}{\E\,Y}
  \rightarrow \prod_{i=1}^k \big(\lambda_i(1+\delta_i)\big)^{m_i}
\]
\end{lemma}

\begin{proof}
As in the previous case 
\[
\E(Y[X_2]_{m_1}\cdots [X_{2k}]_{m_k}) = 
\sum_{S,T} \sum_{C_1,\ldots,C_r} \pr(A_{S,T} \wedge_{i=1}^r A_{C_i}),
\]
where $r=\sum_{i=1}^k m_i$ and in the sum $C_1,\ldots,C_{m_1}$ are different
cycles of length $2$,$C_{m_1+1},\ldots,C_{m_1+m_2}$ are different
cycles of length $4$ etc. 
It is easy to see that the contribution to
the sum coming from the cases where two of the cycles intersect is
$o(1)$. Therefore it suffices to consider disjoint cycles. 

We now repeat the previous argument where
$\eta$ and $\xi$ will refer to $r$ disjoint cycles. 
We then obtain the same formula for $P_2/P_3$, where now $i$ is the
total length of the cycles, $j_1$ is the total number of black
vertices in $V_1$ covered by cycles and $j_2$ is the total number of
black vertices in $V_2$ covered by cycles. 
Finally in order to evaluate the sum corresponding
to~(\ref{eq:gfunc}), we note that it factorizes 
as a power of  the 
sums for individual cycles. 
\end{proof}

\begin{lemma} \label{lem:sumdelta}
\[
\sum_{ {\rm even}\ i\ge 2} \,\lambda_i\,{\delta_i}^2 < \infty
\]
and 
\[
\exp( \sum_{ {\rm even}\ i\ge 2} \,\lambda_i\,{\delta_i}^2 ) = 
\tau^{\a,\b}(d).
\]
\end{lemma}
\begin{proof}
Finding $r(d,i)$ is a well-known problem; one can for example solve
the recurrence 
\[
r(d,i)=d(d-1)^{i-1}-r(d,i-1)
\]
to obtain
$$
r(d,i) = (d-1)^i+(-1)^i(d-1) 
$$
and again we only pay attention to $i$ even. 

We now have
\bean
\sum_{ {\rm even}\ i\ge 2} \,\lambda_i\,{\delta_i}^2
&=& 
\sum_{ {\rm even}\ i\ge 2}\,
\frac{1}{i}\big((d-1)^i+ (d-1) \big)\left(\frac{ \a \b  }{(1-\a) (1-\b) }\right)^{i}\\
&=&\rho\left(\frac{(d-1) \a\b}{(1-\a) (1-\b)}\right)+ (d-1)\rho\left(\frac{  \a\b}{(1-\a) (1-\b)}\right)
\eean
where $\rho(x)=-\frac12(\ln (1-x) +\ln(1+x)) = -\frac12(\ln (1-x^2))$.
We find that
\bean
1-\frac{(d-1)^2 \a^2\b^2}{(1-\a)^2 (1-\b)^2}
&=&
\frac{ (1-\a-\b+d\a\b) (1-\a-\b-(d-2)\a\b)}{(1-\a)^2 (1-\b)^2},\\
1-\frac{ \a^2 \b^2  }{(1-\a) ^2(1-\b)^2 }
&=&
\frac{  (1-\a-\b) (1-\a-\b + 2\a\b)}{(1-\a)^2 (1-\b)^2},
\eean
and hence
$$e^{\sum \lambda_i\,{\delta_i}^2}
 = 
  \tau^{\a,\b}(d)  
$$
as calculated in Theorem~\ref{thm:2nd}.
\end{proof}

\begin{proofof}{Theorem~\ref{thm:YOK}}
Part (i) of the Theorem holds by Lemma~\ref{lem:lam_i}, Part (ii)
holds by Lemma~\ref{lem:deltas} and parts (iii) and (iv) of
Theorem~\ref{sscond} hold by Lemma~\ref{lem:sumdelta}. 
So a.a.s.\ a random graph has independent sets $S$ and
$T$ counted by $Y$.

\end{proofof}

\appendix

\end{document}